\documentclass{amsart}

\usepackage[T1]{fontenc}
\usepackage[leqno]{amsmath}
\usepackage{amssymb, amsthm}
\usepackage{enumerate, latexsym, mathrsfs, mdwlist}

\theoremstyle{definition}
\newtheorem{nul}{}[section]
\newtheorem{dfn}[nul]{Definition}

\newtheorem{ntn}[nul]{Notation}
\newtheorem{exm}[nul]{Example}
\newtheorem{ctrexm}[nul]{Counterexample}

\newtheorem*{ntn*}{Notation}

\theoremstyle{plain}
\newtheorem{thm}[nul]{Theorem}
\newtheorem{prp}[nul]{Proposition}
\newtheorem{lem}[nul]{Lemma}
\newtheorem{cor}{Corollary}[nul]

\numberwithin{equation}{nul}

\DeclareMathOperator{\chr}{char}
\DeclareMathOperator{\colim}{colim}
\DeclareMathOperator{\Ext}{Ext}
\DeclareMathOperator{\id}{id}
\DeclareMathOperator{\Map}{Map}
\DeclareMathOperator{\Mor}{Mor}

\DeclareMathOperator{\rk}{rk}
\DeclareMathOperator{\Shv}{\mathbf{Shv}}
\DeclareMathOperator{\Spec}{Spec}
\DeclareMathOperator{\Sub}{Sub}

\newcommand{\Prod}{\prod}
\newcommand{\Coprod}{\coprod}

\newcommand{\coloneq}{\mathrel{\mathop:}=}

\newcommand{\BB}{\mathbf{B}}
\newcommand{\CC}{\mathbf{C}}

\newcommand{\FF}{\mathbf{F}}

\newcommand{\LL}{\mathbf{L}}

\newcommand{\NN}{\mathbf{N}}

\newcommand{\QQ}{\mathbf{Q}}
\newcommand{\RR}{\mathbf{R}}

\newcommand{\UU}{\mathbf{U}}
\newcommand{\VV}{\mathbf{V}}

\newcommand{\ZZ}{\mathbf{Z}}

\newcommand{\coh}{\mathrm{coh}}
\newcommand{\op}{\mathrm{op}}
\newcommand{\red}{\mathrm{red}}

\usepackage{tikz}
\usetikzlibrary{matrix,arrows,decorations}

\pgfdeclaredecoration{single line}{initial}{\state{initial}[width=\pgfdecoratedpathlength-1sp]{\pgfpathmoveto{\pgfpointorigin}}\state{final}{\pgfpathlineto{\pgfpointorigin}}} 

\newcommand{\fromto}[2]{{#1}\ \tikz[baseline]\draw[>=stealth,->](0,0.5ex)--(0.5,0.5ex);\ {#2}}
\newcommand{\into}[2]{{#1}\ \tikz[baseline]\draw[>=stealth,right hook->](0,0.5ex)--(0.5,0.5ex);\ {#2}}

\newcommand{\goesto}[2]{{#1}\ \tikz[baseline]\draw[|->](0,0.5ex)--(0.5,0.5ex);\ {#2}}

\newcommand{\Alg}{\mathbf{Alg}}
\newcommand{\Cat}{\mathscr{C}\!\mathit{at}}
\newcommand{\et}{\mathrm{\acute{e}t}}
\newcommand{\Et}{\mathrm{\acute{E}t}}
\newcommand{\eet}{\acute{\mathit{e}}\mathit{t}}
\newcommand{\EEt}{\mathbf{\acute{E}t}}
\newcommand{\perf}{\mathrm{perf}}
\newcommand{\Sch}{\mathrm{Sch}}
\newcommand{\Sect}{\mathrm{Sect}}
\newcommand{\Set}{\mathbf{Set}}
\newcommand{\Syn}{\mathrm{Syn}}
\newcommand{\trig}{\mathrm{trig}}
\newcommand{\UH}{\mathrm{UH}}
\newcommand{\wn}{\mathrm{wn}}
\newcommand{\Zar}{\mathrm{Zar}}

\newcommand{\frenchquote}[1]{{\guillemotleft\hspace{0.25em}}{#1}{\hspace{0.25em}\guillemotright}}

\begin{document}

\title{Topological Rigidification of Schemes}

\author{Clark Barwick}
\email{clarkbar@gmail.com}
\address{Department of Mathematics\\
Harvard University\\
One Oxford Street\\
Cambridge, MA 02138\\
USA}

\thanks{This work was supported in part by a Collaborative Research Grant from the National Science Foundation, DMS--0905950.}

\begin{abstract} We show that any of a large class of schemes receives a universal homeomorphism from a reduced scheme that in turn receives no nontrivial universal homeomorphism from any other reduced scheme. This construction serves as a categorical input for the formal inversion of universal homeomorphisms; the result is an $\infty$-category that embeds as a full subcategory of the h $\infty$-topos of Voevodsky.
\end{abstract}

\maketitle

\tableofcontents

\section*{Introduction}

The idea to treat the collection of \emph{universal homeomorphisms} of schemes as a class of \emph{weak equivalences} seems to have been Grothendieck's: a quarter of a century ago, in his \quotedblbase{Brief an Faltings}`` \cite[(8)]{MR1483108}, Grothendieck suggested that, for a field $K$ finitely generated over $\QQ$, after formally inverting universal homeomorphisms, the functor $\goesto{X}{X_{\et}}$ defines a fully faithful embedding of the category of varieties over $K$ into the category of topoi over $(\Spec K)_{\et}$. More generally, Grothendieck expected that a variant of the conjecture should hold in positive characteristic, but only after passage to the perfect closure of the base field. This conjecture was discussed in the \frenchquote{Esquisse d'un program} \cite[Endnote (3), p. 53]{MR1483107}, and, in characteristic zero, a version thereof was proved by Voevodsky \cite{MR1098621}. (See also \cite[\S 3]{MR1931963}.)

The contributions of this paper to this story are: (1) to provide a geometric construction --- the \emph{topological rigidification} --- that yields a model for any quasicompact and quasiseparated scheme up to universal homeomorphism, and (2) to use this construction to invert universal homeomorphisms in the \'etale topos, yielding an $\infty$-topos related to the one corresponding to Voevodsky's h topology. In a sequel to this paper, we aim to prove a generalization of Grothendieck's conjecture, which, in effect, identifies a large class of Artin $n$-stacks of finite type over $\Spec\ZZ$ that can, up to universal homeomorphism, be reconstructed from their \'etale $\infty$-topoi.

We describe in detail the contents of this paper and the manner in which its aims are achieved.

Recall that a morphism $f\colon\fromto{X}{Y}$ of schemes is said to be a \emph{universal homeomorphism} if any base change of $f$ is a homeomorphism on the underlying topological space. We collect in \S \ref{sect:review} a number of elementary facts about this remarkable class of morphisms. The various equivalent characterizations of universal homeomorphisms [\ref{thm:uhequiv}] yield an effective way of testing whether a morphism is a universal homeomorphism. Thinking of universal homeomorphisms as weak equivalences is justified, in part, by the fact that they satisfy the two-out-of-three axiom [\ref{lem:uh2of3}]. Consequently, virtually any question about universal homeomorphisms can be addressed by treating the nilimmersion and schematically dominant cases separately [\ref{lem:uhschdomnilimm}]. We discuss examples and nonexamples of universal homeomorphisms in \ref{exm:univhomeos} and \ref{ctrexm:univhomeos}.

The study of universal homeomorphisms can be reduced to the study of universal homeomorphisms with good finiteness properties. Any universal homeomorphism can be approximated by a finite universal homeomorphism [\ref{prp:cofiltlimuhuh}]. The study of universal homeomorphisms can be further reduced by means of a factorization of any finite universal homeomorphism as a nilimmersion followed by a universal homeomorphism of finite presentation [\ref{prp:factorfiniteuh}]. Universal homeomorphisms of finite presentation can in turn be studied by means of absolute noetherian approximation of the base thanks to \ref{thm:fpresuhsdescend}.

The condition that a morphism be a universal homeomorphism can be regarded as \emph{complementary} to the condition that it be \'etale. The first piece of justification for this manner of thinking comes in \S \ref{sect:propfunet} from the observation that a universal homeomorphism is \'etale if and only if it is an isomorphism [\ref{prp:etuhisiso}].

The second piece of justification for this line of thinking begins with the \frenchquote{\'equivalence remarquable de cat\'egories} of \cite[18.1.2]{MR39:220}: base change along any nilimmersion $\fromto{Z}{X}$ yields an equivalence of categories between the category $\Et_X$ of \'etale $X$-schemes and the category $\Et_Z$ of \'etale $Z$-schemes; thus the morphism of small \'etale topoi $\fromto{Z_{\et}}{X_{\et}}$ is also an equivalence. Grothendieck generalized this result \cite[Exp. IX, 4.10]{MR50:7129} and \cite[Exp. VIII, 1.1]{MR50:7131} to all universal homeomorphisms; see also \cite[5.21]{rydh:submers}. We give a new, less circuitous, proof (without descent theory) of this \frenchquote{invariance topologique} of the \'etale topos in \S \ref{sect:invtopet}.

Frobenius morphisms are a good source of examples of universal homeomorphisms. In fact, in \S \ref{sect:uhFrob}, we follow work of Koll\'ar \cite{MR1432036} to sketch a proof that, in positive characteristic, every finite universal homeomorphism of noetherian schemes is a section up to Frobenius; that is, any universal homeomorphism factors a sufficiently high-powered geometric Frobenius [\ref{prp:kollar}]. Thus, the geometric Frobenii are cofinal in the category of universal homeomorphisms out of a given noetherian scheme.

A scheme is said to be \emph{topologically rigid} if the only universal homeomorphisms to it are retractions of nil-thickenings. In effect, topologically rigid varieties are those varieties with perfect function fields for which any non-normality is as transverse as possible. We introduce this class of objects in \ref{dfn:toprigid}. Filtered limits of topologically rigid schemes remain topologically rigid [\ref{prp:jcommuteswithfiltlims}]; this permits us to employ absolute noetherian reduction arguments in the sequel. The condition of topological rigidity has a purely geometric analogue called \emph{seminormality}; these are equivalent in characteristic zero, but in positive characteristic, they differ dramatically, as our examples and nonexamples show [\ref{exm:twnsch}, \ref{ctrexm:twnsch}].

The inclusion functor from the category of reduced schemes into the category of all schemes admits a right adjoint $\goesto{Y}{Y_{\red}}$. Similarly, in \ref{prp:trigexists} we show that the inclusion functor from the category of topologically rigid schemes into the category of all schemes admits a right adjoint $\goesto{Y}{Y_{\trig}}$; this construction is called \emph{topological rigidification}. As $Y_{\red}$ can be thought of as the best reduced approximation to $Y$, so $Y_{\trig}$ can be thought of as the best topologically rigid approximation to $Y$.

It is a remarkable consequence of Koll\'ar's result that the topological rigidification of a reduced noetherian $\FF_p$-scheme coincides with its perfection; in effect, the limit over all universal homeomorphisms can be replaced by the system of Frobenii [\ref{exm:trigvsperf}]. This computation of the topological rigidification illustrates two peculiar phenomena in positive characteristic. The first is the ability the Frobenius has to resolve problems of transversality. The second is the price one has to pay for this simplicity: the topological rigidification of a scheme with pleasant finiteness properties is unlikely to retain those good properties. It is for this reason that we are forced to contend with such general categories of schemes. (We note that Voevodsky alludes to the existence of $\goesto{Y}{Y_{\trig}}$ after \cite[3.2.10]{MR1403354}, but he does not investigate the construction further, precisely because of its failure to preserve finiteness properties.)

In \S \ref{sect:wknorm}, following the early work of Andreotti--Bombieri \cite{MR0266923} and the later work of Manaresi \cite{MR556308}, we sketch proofs of the existence and properties of a relative version of the topological rigidification --- the \emph{weak normalization} of a scheme $Y$ relative to a morphism $\fromto{Z}{Y}$. This yields a powerful description of the topological rigidification of reduced schemes with locally finitely many irreducible components [\ref{cor:computetrig}], and it permits us to prove \ref{cor:etcovtrigistrig}, which is key for the construction of the final section.

At last, in \S \ref{sect:htopos}, we invert the universal homeomorphisms in the \'etale $\infty$-topos. Our method of doing so makes essential use of the topological rigidification functor. The result is called the \emph{topologically invariant étale $\infty$-topos}, and it is closely related to the h topology of Voevodsky \cite{MR1403354}. Topologically rigid schemes thus form a generating site for the h $\infty$-topos for which the topology is subcanonical. We conclude with a natural isomorphism between the nonabelian $h$-cohomology of a group scheme $G$ and the nonabelian \'etale cohomology of $G_{\trig}$ [\ref{cor:hcohomology}], giving a new proof of a theorem of Suslin--Voevodsky.

\subsection*{Acknowledgements} Thanks to Dustin Clausen and Kirsten Wickelgren for questions, comments, and corrections about this work.

\subsection*{Conventions} The language of $\infty$-categories and $\infty$-topoi \cite{lurie_inftytopoi} is restricted to the final section. There, the assertions have been formulated so that one may simply remove all instances of the phrases ``simplicial'' and ``up to homotopy'' as well as the symbol ``$\infty$'' in order to obtain correct statements and proofs about the (ordinary) h topos; however, in the sequel to this paper, the use of $\infty$-topoi will not be so readily avoidable.

We use the Universe Axiom of \cite[Exp I, \S 0]{MR50:7130}, and we fix two universes $\UU\in\VV$. All rings, modules, schemes, etc., will be $\UU$-small. We employ the following notational conventions for categories or $\infty$-categories:
\begin{enumerate}[{(N.}1{)}]
\item Roman characters $A$, $B$, $C$, \dots, etc., will denote categories and $\infty$-categories that are essentially $\VV$-small and locally $\UU$-small.
\item Bold characters $\mathbf{A}$, $\BB$, $\CC$, \dots, etc., will denote categories and $\infty$-categories that are locally $\VV$-small.
\end{enumerate}

\section{Review: Universal homeomorphisms}\label{sect:review}

\setcounter{nul}{-1}

\begin{ntn} For a scheme $X$, denote by $|X|$ the underlying topological space; for a morphism $f\colon\fromto{X}{Y}$ of schemes, denote by $|f|\colon\fromto{|X|}{|Y|}$ the underlying continuous map. We say that $f$ is \emph{surjective}, \emph{injective}, \emph{closed}, \emph{open}, etc., if and only if $|f|$ is so. Denote by $f^{\sharp}\colon\fromto{\mathscr{O}_Y}{f_{\star}\mathscr{O}_X}$ the corresponding morphism of sheaves.
\end{ntn}

\begin{thm}[\protect{\cite[18.12.8]{MR39:220}, \cite[3.8.10]{EGA1}}]\label{thm:entierequiv} The following conditions on a morphism $f\colon\fromto{Y}{X}$ of schemes are equivalent.
\begin{enumerate}[(\ref{thm:entierequiv}.1)]
\item The morphism $f$ is integral \emph{(entier)}.
\item The morphism $f$ is affine and the $\mathscr{O}_X$-algebra $f_{\star}\mathscr{O}_Y$ is integral.
\item The morphism $f$ is affine and universally closed.
\item The morphism $f$ is affine, and for any set $S$, the induced morphism
\begin{equation*}
\fromto{\mathbf{A}^{\!S}_Y\coloneq\Spec_Y(\mathscr{O}_Y[S])}{\Spec_X(\mathscr{O}_X[S])=:\mathbf{A}^{\!S}_X}
\end{equation*}
is closed.
\end{enumerate}
\end{thm}

\begin{lem}[\protect{\cite[6.1.5(v)]{MR36:177b}}]\label{lem:gfintonlyiffint} Suppose $\fromto{Z}{X}$ an integral morphism of schemes, and suppose $\fromto{Y}{X}$ a separated morphism of schemes. Then any morphism $\fromto{Z}{Y}$ over $X$ is integral as well.
\end{lem}

\begin{dfn}\label{dfn:seprank} Suppose $f\colon\fromto{Y}{X}$ a quasifinite morphism of schemes, and suppose $x\in X$ a point. Following \cite[I, 6.5.9]{EGA1}, we define the \emph{separable rank} of $Y$ over $x$ to be the natural number
\begin{equation*}
\rk_x(Y)\coloneq\sum_{y\in f^{-1}(x)}[\kappa(y):\kappa(x)]_{\mathrm{s}}.
\end{equation*}
\end{dfn}

\begin{thm}[\protect{\cite[3.7.1]{EGA1}}]\label{thm:uiequiv} The following conditions on a morphism $f\colon\fromto{Y}{X}$ of schemes are equivalent.
\begin{enumerate}[(\ref{thm:uiequiv}.1)]
\item The morphism $f$ is universally injective \emph{(radiciel)}.
\item The map $|f|$ is injective, and for any point $y\in Y$, if $x=f(y)$, the residue field extension $\into{\kappa(x)}{\kappa(y)}$ is purely inseparable.
\item For any field $k$, the induced map $\fromto{Y(k)}{X(k)}$ is injective.
\item For any field $k$, there is an algebraically closed extension $\Omega$ of $k$ such that the induced map $\fromto{Y(\Omega)}{X(\Omega)}$ is injective.
\item The diagonal $|\Delta_{Y/X}|\colon\fromto{|Y|}{|Y\times_XY|}$ is surjective.
\suspend{enumerate}
If, in addition, $f$ is quasifinite, then the following condition may be added.
\resume{enumerate}[{[(\ref{thm:uiequiv}.1)]}]
\item For any $x\in|X|$, the separable rank $\rk_x(Y)=1$.
\end{enumerate}
\end{thm}

\begin{cor}\label{cor:gfuionlyiffui} Suppose $\fromto{Z}{X}$ a universally injective morphism of schemes. Then any morphism $\fromto{Z}{Y}$ of schemes over $X$ is universally injective as well.
\end{cor}

\begin{thm}[\protect{\cite[18.12.11]{MR39:220}}]\label{thm:uhequiv} The following conditions on a morphism $f\colon\fromto{Y}{X}$ of schemes are equivalent.
\begin{enumerate}[(\ref{thm:uhequiv}.1)]
\item The morphism $f$ is a universal homeomorphism.
\item The morphism $f$ is surjective, universally injective, and universally closed.
\item The morphism $f$ is surjective, universally injective, and integral.
\item\label{item:uhexplicit} The morphism $f$ satisfies the following pair of conditions.
\begin{enumerate}[(\ref{thm:uhequiv}.\ref{item:uhexplicit}.1)]
\item For any $x\in|X|$, the fiber $f^{-1}(x)=\{y\}$ such that the residue field extension $\into{\kappa(x)}{\kappa(y)}$ is purely inseparable.
\item For any set $S$, the induced morphism $\fromto{\mathbf{A}^{\!S}_Y}{\mathbf{A}^{\!S}_X}$ is closed.
\end{enumerate}
\end{enumerate}
\end{thm}

\begin{prp}\label{lem:uh2of3} Suppose $X$ a scheme. The set of universal homeomorphisms over $X$ satisfies the two-out-of-three axiom. More precisely, the following obtains.
\begin{enumerate}[(\ref{lem:uh2of3}.1)]
\item\label{item:gfuhonlyiffuh} If $\fromto{W}{U}$ is a universal homeomorphism of $X$-schemes and $\fromto{V}{U}$ is a separated, surjective, and universally injective morphisms of schemes, then any morphism $\fromto{W}{V}$ over $U$ is a universal homeomorphism as well.
\item If $\fromto{W}{V}$ a universal homeomorphism of $X$-schemes, then a morphism $\fromto{V}{U}$ of $X$-schemes is a universal homeomorphism if and only if the composite $\fromto{W}{U}$ is so.
\end{enumerate}
\begin{proof} (\ref{lem:uh2of3}.\ref{item:gfuhonlyiffuh}) follows immediately from \ref{cor:gfuionlyiffui} and \ref{lem:gfintonlyiffint}.

The second claim is the observation that for any morphism $\fromto{U'}{U}$ of $X$-schemes, then $\fromto{V\times_UU'}{U'}$ is a homeomorphism if and only if the composite $\fromto{W\times_UU'}{U'}$ is so.
\end{proof}
\end{prp}

\begin{cor}\label{lem:uhschdomnilimm} The class of universal homeomorphisms $f\colon\fromto{Y}{Y'}$ admits functorial factorizations into schematically dominant universal homeomorphisms $\fromto{Y}{Y''}$ followed by nilimmersions $\into{Y''}{Y'}$.
\begin{proof} Simply factor $f$ through its schematic image \cite[I 6.10.1]{EGA1}.
\end{proof}
\end{cor}

\begin{prp}\label{prp:schdomuhsareepis} A schematically dominant universal homeomorphism is an epimorphism in the category of schemes.
\begin{proof} Any such morphism $f\colon\fromto{X}{Y}$ induces a homeomorphism $|f|$ of topological spaces and a monomorphism $f^{\sharp}\colon\fromto{\mathscr{O}_Y}{f_{\star}\mathscr{O}_X}$ of sheaves.
\end{proof}
\end{prp}

\begin{exm}\label{exm:univhomeos}
\begin{enumerate}[(\ref{exm:univhomeos}.1)]
\item Any nilimmersion is a universal homeomorphism, and any section of a universal homeomorphism is a nilimmersion.
\item If $\into{K}{L}$ is a purely inseparable field extension, then the morphism $\fromto{\Spec L}{\Spec K}$ is a universal homeomorphism. Moreover, for any $K$-scheme $X$, the morphism $\fromto{X\otimes_KL}{X}$ is a universal homeomorphism.
\item Suppose $X$ a reduced scheme whose set of irreducible components is locally finite.\footnote{This means that every quasicompact open of $X$ contains only finitely many irreducible components. This condition holds, e.g., if the topological space $|X|$ is locally noetherian.} Denote by $X^{\nu}$ its normalization. Then the following are equivalent: (1) the natural morphism $\fromto{X^{\nu}}{X}$ is a universal homeomorphism; (2) $X$ is geometrically unibranch \cite[6.15.1]{MR33:7330}; (3) for every point $x\in X$, the strict henselianization ${}^{\mathrm{sh}}\mathscr{O}_{X,x}$ of the local ring at $x$ is irreducible \cite[18.8.15(c)]{MR39:220}.
\item Suppose $k$ an algebraically closed field of characteristic $2$. The subring $k[x^2,xy,y]\subset k[x,y]$ induces a morphism
\begin{equation*}
f\colon\fromto{\Spec k[x,y]}{\Spec k[x^2,xy,y]}
\end{equation*}
that is finite, surjective, and an isomorphism away from the $x$-axis. One sees that the ideal $(y)$ lies over the ideal $(xy,y)$; the induced morphism of subschemes
\begin{equation*}
\fromto{\Spec k[x]}{\Spec k[x^2]}
\end{equation*}
is a bijection inducing an isomorphism of residue fields on every closed point. On the generic point of this subscheme, the field extension is $k(x^2)\subset k(x)$. One concludes that $\fromto{Y}{X}$ is a universal homeomorphism.
\end{enumerate}
\end{exm}

\begin{ctrexm}\label{ctrexm:univhomeos}\begin{enumerate}[(\ref{ctrexm:univhomeos}.1)]
\item Consider the subring $\RR[x,\sqrt{-1}x]\subset\CC[x]$. The induced morphism $f\colon\fromto{\Spec\CC[x]}{\Spec\RR[x,\sqrt{-1}x]}$ is a bijection, and away from the origin it is an isomorphism. At the origin, the residue field extension is $\RR\subset\CC$; hence $f$ is not a universal homeomorphism.
\item Suppose $k$ an algebraically closed field of characteristic $0$. Suppose $X$ a nodal curve over $k$, and let $\fromto{X^{\nu}}{X}$ be its normalization, so that over the singular point $x\in X(k)$ there are two points $u,v\in X^{\nu}(k)$. Remove $u$, and consider the restriction $\fromto{X^{\nu}-\{u\}}{X}$; this is a bijection on all points that induces an isomorphism on all residue field extensions. It is neither a homeomorphism, nor an integral morphism.
\end{enumerate}
\end{ctrexm}

\section{Limits of universal homeomorphsims and approximation}\label{sect:limsapprox}

\begin{prp}\label{prp:cofiltlimuhuh} Suppose $X$ a scheme. Suppose $\Lambda$ a small filtered category, and suppose $W_{\ast}\colon\fromto{\Lambda^{\op}}{(\Sch/X)}$ a diagram of $X$-schemes such that for any object $\alpha\in\Lambda^{\op}$, the structure morphism $p_{\alpha}\colon\fromto{W_{\alpha}}{X}$ is a universal homeomorphism. Then the natural morphism
\begin{equation*}
p\colon\fromto{W\coloneq\lim_{\alpha\in\Lambda^{\op}}W_{\alpha}}{X}
\end{equation*}
is a universal homeomorphism.
\begin{proof} All the bonding morphisms $\fromto{W_{\alpha}}{W_{\beta}}$ are universal homeomorphisms by \ref{lem:uh2of3}. It follows from \cite[8.3.8(i)]{MR36:178} that $p$ is surjective. For any field $k$, the diagram $W(k)\colon\fromto{\Lambda^{\op}}{\Set}$ is a diagram of injections, whence for any $\alpha\in\Lambda^{\op}$, the map $\fromto{W(k)}{W_{\alpha}(k)}$ is an injection; thus $p$ is a universal injection. It thus remains to show that $p$ is integral. Since $W$ is a diagram of affine $X$-schemes, it is enough to observe that the filtered colimit $\colim_{\alpha\in\Lambda}p_{\alpha,\star}\mathscr{O}_{W_{\alpha}}$ is an integral $\mathscr{O}_X$-algebra.
\end{proof}
\end{prp}

\begin{lem}\label{prp:uhlimcofiltptype} Suppose $X$ a coherent (i.e., quasicompact and quasiseparated) scheme, and suppose $p\colon\fromto{W}{X}$ a universal homeomorphism. Then there exists a filtered category $\Lambda$ and a diagram $W_{\ast}\colon\fromto{\Lambda^{\op}}{(W/\Sch/X)}$ of $X$-schemes under $W$ satisfying the following conditions.
\begin{enumerate}[(\ref{prp:uhlimcofiltptype}.1)]
\item The natural morphism $\fromto{W}{\lim_{\alpha\in\Lambda^{\op}}W_{\alpha}}$ is an isomorphism.
\item For any $\alpha\in\Lambda^{\op}$, the morphism $p_{\alpha}\colon\fromto{W_{\alpha}}{X}$ is a finite universal homeomorphism.
\end{enumerate}
\begin{proof} Write $\mathscr{A}$ for the quasicoherent $\mathscr{O}_X$-algebra $p_{\star}\mathscr{O}_W$, so that $W=\Spec_X\mathscr{A}$. By \cite[I 6.9.15]{EGA1}, write $\mathscr{A}$ as a filtered colimit $\colim_{\alpha\in\Lambda}\mathscr{A}_{\alpha}$ of its sub-$\mathscr{O}_X$-algebras of finite type. For any $\alpha\in\Lambda$, set $W_{\alpha}\coloneq\Spec_X\mathscr{A}_{\alpha}$. The morphism $\fromto{W}{W_{\alpha}}$ is integral and universally injective; moreover, it is schematically dominant, hence surjective. Hence it is a universal homeomorphism, and by \ref{lem:uh2of3}, so is $p_{\alpha}$.
\end{proof}
\end{lem}

\begin{prp}\label{prp:factorfiniteuh} Suppose $X$ a coherent scheme, and suppose $p\colon\fromto{W}{X}$ any finite universal homeomorphism. Then there exists a universal homeomorphism $p'\colon\fromto{W'}{X}$  of finite presentation and a nilimmersion $\into{W}{W'}$ over $X$.
\begin{proof} By \cite[4.2]{MR2356346}, $p$ factors as a separated morphism $\fromto{Z}{X}$ of finite presentation followed by a closed immersion $i\colon\into{W}{Z}$ over $X$. Set $\mathscr{I}\coloneq\ker[\fromto{\mathscr{O}_Z}{i_{\star}\mathscr{O}_W}]$. Let $\Lambda$ denote the category of quasicoherent subideals $\mathscr{J}_{\alpha}\subset\mathscr{I}$ of finite type. Now \cite[I 6.9.16(iii)]{EGA1} $\Lambda$ is filtered, and the $\mathscr{O}_Z$-algebra $\mathscr{O}_Z/\mathscr{J}_{\alpha}$ is of finite presentation for any object $\alpha\in\Lambda$. Write $W_{\alpha}:=\Spec_Z(\mathscr{O}_Z/\mathscr{J}_{\alpha})$; hence $\lim_{\alpha\in\Lambda}W_{\alpha}\cong W$. For any $\alpha\in\Lambda$, $\fromto{W_{\alpha}}{X}$ is surjective, $\into{W}{W_{\alpha}}$ is a closed embedding, and $\into{W_{\alpha}}{Z}$ is a closed embedding. It remains to show that for some $\alpha_0\in\Lambda$, the morphism $p_{\alpha_0}$ is finite and universally injective. Suppose $\mathscr{V}$ a finite affine open covering of $X$.

(1) First, for every element $V\in\mathscr{V}$, the scheme $V\times_XW$ is affine; hence there is an $\alpha_{-2}\in\Lambda$ such that for any $V\in\mathscr{V}$, the scheme $V\times_XW_{\alpha_{-2}}$ is affine. Thus for every $\alpha\in(\alpha_{-2}/\Lambda)$, the morphism $\fromto{W_{\alpha}}{X}$ is affine.

(2) Next, for every $V\in\mathscr{V}$, the $\Gamma(V,\mathscr{O}_X)$-algebra $\Gamma(V\times_XW_{\alpha_{-2}},\mathscr{O}_Z)$ is finitely presented; choose a finite set $S_V$ of generators thereof. The images of these generators in $\Gamma(V\times_XW_{\alpha},\mathscr{O}_Z)$ are generators for any $\alpha\in(\alpha_{-2}/\Lambda)$, and since $p$ is finite, the image of each generator $x\in S_V$ in $\Gamma(V\times_XW,\mathscr{O}_Z)$ satisfies a monic polynomial $f_{V,x}\in\Gamma(V,\mathscr{O}_X)[t]$. Hence there exists an $\alpha_{-1}\in(\alpha_{-2}/\Lambda)$ such that the image of $f(V,x)(x)$ vanishes in $\Gamma(V\times_XW_{\alpha_{-1}},\mathscr{O}_Z)$ for all $V\in\mathscr{V}$ and all $x\in S_V$. Thus $p_{\alpha_{-1}}$ is finite, and for every $\alpha\in(\alpha_{-1}/\Lambda)$, the morphism $p_{\alpha}\colon\fromto{W_{\alpha}}{X}$ is finite.

(3) Finally, for any $\alpha\in(\alpha_{-1}/\Lambda)$, set
\begin{equation*}
Z_{\alpha}\coloneq\{z\in Z\ |\ \rk_z(W_{\alpha})=1\}, 
\end{equation*}
where $\rk_z(W_{\alpha})$ denotes the separable rank [\ref{dfn:seprank}] of $W_{\alpha}$ over $z$. By \cite[9.7.9]{MR36:178}, the sets $Z_{\alpha}$ are constructible subsets of $Z$. One verifies that
\begin{equation*}
Z=\bigcup_{\alpha\in(\alpha_{-1}/\Lambda)}Z_{\alpha},
\end{equation*}
a filtered union. Since $Z$ is quasicompact, there is an $\alpha_0\in\Lambda$ such that $Z=Z_{\alpha_0}$, so for any $\alpha\in(\alpha_0/\Lambda)$, the surjection $\fromto{W_{\alpha}}{Z}$ is universally injective, and for every $\alpha\in(\alpha_0/\Lambda)$, the morphism $p_{\alpha}\colon\fromto{W_{\alpha}}{X}$ is universally injective.
\end{proof}
\end{prp}

\begin{ntn*} For the next result, suppose $X_0$ a scheme, suppose $\Lambda$ a filtered category, and suppose $\mathscr{A}\colon\fromto{\Lambda}{\Alg(\mathscr{O}_{X_0})}$ a diagram of quasicoherent commutative $\mathscr{O}_{X_0}$-algebras. For any $\alpha\in\Lambda$, write $X_{\alpha}:=\Spec_{X_0}\mathscr{A}_{\alpha}$. Suppose $Y_0$ and $Z_0$ two $X_0$-schemes; for any $\alpha\in\Lambda$, write $Y_{\alpha}:=Y_0\times_{X_0}X_{\alpha}$ and $Z_{\alpha}:=Z_0\times_{X_0}X_{\alpha}$. Finally, write $X:=\lim_{\alpha\in\Lambda^{\op}}X_{\alpha}$, $Y:=\lim_{\alpha\in\Lambda^{\op}}Y_{\alpha}$, and $Z:=\lim_{\alpha\in\Lambda^{\op}}Z_{\alpha}$.
\end{ntn*}

\begin{thm}\label{thm:fpresuhsdescend} Suppose $X_0$ quasicompact, and suppose $Y_0$ and $Z_0$ locally of finite presentation over $X_0$. Then any universal homeomorphism $\fromto{Y}{Z}$ over $X$ descends to a universal homeomorphism $\fromto{Y_{\alpha}}{Z_{\alpha}}$ over $X_{\alpha}$ for some $\alpha\in\Lambda$.
\begin{proof} Apply \cite[8.10.5(vi, vii, x)]{MR36:178}.
\end{proof}
\end{thm}

\section{The \frenchquote{propri\'et\'e fondamentale} of \'etale morphisms}\label{sect:propfunet}

\begin{prp}\label{prp:etuhisiso} A morphism $f\colon\fromto{Y}{X}$ of schemes is an \'etale universal homeomorphism if and only if it is an isomorphism.
\begin{proof} The sufficiency is obvious. For the necessity, suppose $f$ an \'etale universal homeomorphism. Then any section of $f$ will be a surjective open immersion \cite[17.4.1(b${}^{\prime\prime}$)]{MR39:220}, hence an isomorphism. It suffices to construct such a section after a faithfully flat base change. So consider the projection $\fromto{Y\times_XY}{Y}$; the diagonal $\Delta_{Y/X}\colon\fromto{Y}{Y\times_XY}$ is the desired section.
\end{proof}
\end{prp}

\begin{cor}[\protect{\cite[17.9.1]{MR39:220}}] A morphism $f\colon\fromto{Y}{X}$ of schemes is \'etale and universally injective if and only if it is an open immersion.
\begin{proof} The sufficiency is obvious. So suppose $f$ \'etale and universally injective. Since $f$ is flat and locally of finite presentation, it is universally open \cite[2.4.6]{MR33:7330}. It is thus a universal homeomorphism onto its image; one now applies the proposition.
\end{proof}
\end{cor}

\begin{cor} A morphism $f\colon\fromto{Y}{X}$ of schemes is unramified and universally injective if and only if it is a monomorphism.
\begin{proof} Factor $f$ through its schematic image, and apply the previous corollary.
\end{proof}
\end{cor}

\begin{cor}[\protect{\cite[17.9.3]{MR39:220}}]\label{cor:sectsareopensisoover} Suppose $f\colon\fromto{Y}{X}$ an \'etale (respectively, \'etale and separated) morphism of schemes. Then every section of $f$ is an open (resp. open and closed) immersion, and the set $\Sect_X(Y)$ of sections of $f$ is in canonical bijection with the set of open subsets (resp. connected components) $U\subset|Y|$ such that $f|_U\colon\fromto{U}{X}$ is a universal homeomorphism and hence an isomorphism.
\end{cor}

\section{The \frenchquote{invariance topologique} of the \'etale topos}\label{sect:invtopet}

\begin{dfn}\label{dfn:Eequivalence} Suppose $E$ a presheaf of categories on the category $\Sch$ of schemes.\footnote{By this we mean a pseudofunctor $\fromto{\Sch^{\op}}{\Cat_{\VV}}$ valued in $\VV$-small $\UU$-categories or, equivalently, a cartesian fibration (with fibers in categories) over $\Sch$.} Then a morphism $\fromto{Y}{X}$ of schemes will be called an \emph{$E$-equivalence} if the induced functor $\fromto{E_X}{E_Y}$ is an equivalence of categories; it is called a \emph{universal $E$-equivalence} if, for any morphism $\fromto{X'}{X}$, the induced functor $\fromto{E_{X'}}{E_{X'\times_XY}}$ is an equivalence.
\end{dfn}

\begin{ntn} For any scheme $X$, write $\Et_X$ for the category of \'etale $X$-schemes; then pullback defines a presheaf of categories $\Et$ on the category $\Sch$ of schemes. Write $\tau_{\leq 1}X_{\et}$ for the ordinary small \'etale topos of $X$ --- i.e., the category of sheaves of sets on $\Et_X$; then pullback defines a presheaf of categories $\eet$ on the category $\Sch$ of schemes.
\end{ntn}

\begin{prp}\label{prp:etaletopoigiveZar} For any schemes $X$ and $Y$ and any morphism $\phi\colon\fromto{\tau_{\leq 1}X_{\et}}{\tau_{\leq 1}Y_{\et}}$ of small \'etale topoi, there is a unique continuous map $f\colon\fromto{|X|}{|Y|}$ such that for any open set $U\subset|Y|$, one has $f^{-1}U=\phi^{\star}U$.
\begin{proof} The underlying space of any scheme is sober, so it is enough to observe that such a morphism $\phi$ induces a homomorphism
\begin{equation*}
\fromto{\Zar_Y\cong\Sub(\tau_{\leq 1}Y_{\et})}{\Sub(\tau_{\leq 1}X_{\et})\cong\Zar_X}
\end{equation*}
of Zariski locales.
\end{proof}
\end{prp}

\begin{cor}\label{lem:Etequivsarehomeos} An $\eet$-equivalence is a homeomorphism.
\begin{proof} Apply \ref{prp:etaletopoigiveZar} to obtain the continuous inverse.
\end{proof}
\end{cor}

\begin{prp}\label{prp:uhinduceeqetale} For any universal homeomorphism $f\colon\fromto{Y}{X}$ of schemes, the induced functor $\fromto{\Et_X}{\Et_Y}$ is fully faithful.
\begin{proof} By forming pullbacks, it suffices to show that for any universal homeomorphism $\fromto{Y}{X}$ and any \'etale $X$-scheme $X'$, the natural map
\begin{equation*}
\fromto{\Sect_X(X')}{\Sect_Y(Y')},
\end{equation*} 
is a bijection, where $Y'\coloneq X'\times_XY$. Now one may apply \ref{cor:sectsareopensisoover}: since the projection $\fromto{X'\times_XY}{X'}$ is a universal homeomorphism, the set of open subsets $U\subset X'$ such that the induced morphism $\fromto{U}{X}$ is a universal homeomorphism is in canonical bijection with the set of open subsets $V\subset|Y'|$ such that the induced morphism $\fromto{V}{Y}$ is a universal homeomorphism.
\end{proof}
\end{prp}

\begin{cor} \'Etale morphisms satisfy the unique right lifting property with respect to universal homeomorphisms. That is, for any commutative square
\begin{equation*}
\begin{tikzpicture} 
\matrix(m)[matrix of math nodes, 
row sep=4ex, column sep=4ex, 
text height=1.5ex, text depth=0.25ex] 
{Z&Y\\ 
W&X\\}; 
\path[>=stealth,->,font=\scriptsize] 
(m-1-1) edge (m-1-2) 
edge (m-2-1) 
(m-1-2) edge (m-2-2) 
(m-2-1) edge (m-2-2); 
\end{tikzpicture}
\end{equation*}
of schemes in which $\fromto{Y}{X}$ is \'etale and $\fromto{Z}{W}$ is a universal homeomorphism, there exists a unique lift $\fromto{W}{Y}$.
\begin{proof} Pulling back along $\fromto{W}{X}$, it suffices to consider the case in which $W=X$. Hence the claim is that for any scheme $X$ and any \'etale $X$-schemes $Y$ and $Z$, the canonical map
\begin{equation*}
\fromto{\Mor_X(X,Y)}{\Mor_X(Z,Y)}
\end{equation*}
is a bijection. By \ref{prp:uhinduceeqetale}, this amounts to showing that the canonical map
\begin{equation*}
\phi\colon\fromto{\Mor_Z(Z,Y\times_XZ)}{\Mor_Z(Z\times_XZ,Y\times_XZ)}
\end{equation*}
is a bijection. Composition with the diagonal $\Delta_{Z/X}\colon\fromto{Z}{Z\times_XZ}$ defines an inverse $\Delta_{Z/X}^{\star}$ to $\phi$, and $\Delta_{Z/X}^{\star}$ is a bijection since $\Delta_{Z/X}$ is a nilimmersion.
\end{proof}
\end{cor}

\begin{cor}\label{cor:etequivloconbas} The question of whether an \'etale $Y$-scheme descends along a universal homeomorphism $\fromto{Y}{X}$ to an \'etale $X$-scheme is local on $X$.
\begin{proof} Suppose $Z$ an \'etale $Y$-scheme, and suppose $\{U_j\}_{j\in J}$ an affine open cover of $X$; write $\{V_j\}_{j\in J}$ for the pulled back affine open cover of $Y$. Suppose further that for every $j\in J$, one has an \'etale $U_j$-scheme $U'_j$ with an isomorphism
\begin{equation*}
\phi_j\colon V'_j\coloneq V_j\times_YY'\cong V_j\times_{U_j}U'_j.
\end{equation*}
By \ref{prp:uhinduceeqetale}, for any $j,k\in J$, there is a unique isomorphism $U'_j\times_{U_j}U_{jk}\cong U'_k\times_{U_k}U_{jk}$ corresponding to the canonical isomorphism $V'_j\times_{V_j}V_{jk}\cong V'_k\times_{V_k}V_{jk}$. These isomorphisms clearly satisfy the cocycle condition, so gluing the $U'_j$ together along these isomorphisms, one obtains an \'etale $X$-scheme $X'$ with the property that the $\phi_j$ glue together to give an isomorphism $Y'\cong X'\times_XY$.
\end{proof}
\end{cor}

\begin{cor}\label{cor:limitarguhsareets} Suppose $X$ a coherent scheme. Suppose $\Lambda$ a filtered category, and suppose $W\colon\fromto{\Lambda^{\op}}{(\Sch/X)}$ a diagram of $X$-schemes such that for any object $\alpha\in\Lambda^{\op}$, the structure morphism $p_{\alpha}\colon\fromto{W_{\alpha}}{X}$ is a universal homeomorphism that induces an essentially surjective functor $\fromto{\Et_X}{\Et_{W_{\alpha}}}$. Then the natural morphism
\begin{equation*}
\fromto{W\coloneq\lim_{\alpha\in\Lambda^{\op}}W_{\alpha}}{X}
\end{equation*}
is an $\Et$-equivalence.
\begin{proof} The functor $\fromto{\Et_X}{\Et_W}$ is fully faithful by \ref{prp:cofiltlimuhuh} and \ref{prp:uhinduceeqetale}. Now for any \'etale morphism $\fromto{Y}{W}$, apply \cite[8.8.2]{MR36:178} and \cite[17.7.8]{MR39:220} to descend $Y$ to an \'etale $W_{\alpha}$-scheme $Y_{\alpha}$ for some object $\alpha\in\Lambda$; hence the functor $\fromto{\Et_{W_{\alpha}}}{\Et_W}$ --- and thus the functor $\fromto{\Et_X}{\Et_W}$ --- is essentially surjective.
\end{proof}
\end{cor}

\begin{cor}\label{cor:checkonlocals} A universal homeomorphism $f\colon\fromto{Y}{X}$ is an $\Et$-equivalence if, for any point $x\in X$, the universal homeomorphism $\fromto{Y\times_X\Spec\mathscr{O}_{X,x}}{\Spec\mathscr{O}_{X,x}}$ induces a fully faithful functor $\fromto{\Et_{\Spec\mathscr{O}_{X,x}}}{\Et_{Y\times_X\Spec\mathscr{O}_{X,x}}}$.
\begin{proof} The functor $\fromto{\Et_X}{\Et_Y}$ is fully faithful by \ref{prp:uhinduceeqetale}. Write $\Spec\mathscr{O}_{X,x}$ as the cofiltered limit of open affines $\{V_{\alpha}\}_{\alpha\in\Lambda}$ containing $x$. Now for any \'etale $Y$-scheme $Z$, the \'etale $Y\times_X\Spec\mathscr{O}_{X,x}$-scheme $Z\times_X\Spec\mathscr{O}_{X,x}$ descends to an \'etale $\mathscr{O}_{X,x}$-scheme $W_x$. Now apply \cite[8.8.2]{MR36:178} and \cite[17.7.8]{MR39:220} to descend $W_x$ to an \'etale $V_{\alpha}$-scheme $W_{\alpha}$ for some object $\alpha\in\Lambda$. Executing this procedure for every point $x\in X$, one obtains an affine open cover $\{U_{\alpha}\}_{\alpha\in A}$ of $X$ along with an \'etale $U_{\alpha}$-scheme $W_{\alpha}$ and an isomorphism $W_{\alpha}\times_XY\cong Z\times_XU_{\alpha}$. It now follows from \ref{cor:etequivloconbas} that there is an \'etale $X$-scheme $W$ along with an isomorphism $W\times_XY\cong Z$. Thus the functor $\fromto{\Et_X}{\Et_Y}$ is essentially surjective.
\end{proof}
\end{cor}

\begin{prp}\label{prp:nilEtsurj} Nilimmersions are $\Et$-equivalences.
\begin{proof} In view of \ref{prp:uhinduceeqetale}, our claim is that for any nilimmersion $f\colon\fromto{Y}{X}$ of schemes, the induced functor $\fromto{\Et_X}{\Et_Y}$ is essentially surjective. By \ref{cor:etequivloconbas}, one reduces to the case of $X=\Spec A$ and $Y=\Spec B$. Write $I=\ker[\fromto{A}{B}]$ for the nil ideal defining $Y$. Now write $A$ as the filtered colimit $\colim_{\alpha\in\Lambda}A_{\alpha}$ of its subrings of finite type over $\ZZ$, and set $X_{\alpha}\coloneq\Spec A_{\alpha}$. Now set
\begin{equation*}
B_{\alpha}\coloneq A_{\alpha}/(I\cap A_{\alpha})\textrm{\quad and\quad}Y_{\alpha}\coloneq\Spec B_{\alpha}.
\end{equation*}
Thus $Y\cong\lim_{\alpha\in\Lambda^{\op}}Y_{\alpha}$. For any \'etale morphism $p\colon\fromto{Y'}{Y}$, apply \cite[8.8.2]{MR36:178} and \cite[17.7.8]{MR39:220} to descend $Y'$ to an \'etale $Y_{\alpha}$-scheme $Y'_{\alpha}$ for some object $\alpha\in\Lambda$.

One is thus reduced to the case in which $X$, $Y$, and $Y'$ are all noetherian. In this case, $\fromto{Y}{X}$ is nilpotent; by induction, one may further assume that it is square zero. The ideal $I$ can then be regarded as a module on $Y$, and by \cite[III 2.1.7]{MR58:10886a}, there is an obstruction
\begin{equation*}
\varpi(Y'/Y,X)\in\Ext^2(\LL_{Y'/Y},p^{\star}I)
\end{equation*}
whose vanishing is necessary and sufficient for the existence of a deformation of $p$. As $\LL_{Y'/Y}\simeq 0$, we are done.
\end{proof}
\end{prp}

\begin{nul} If $\Syn$ denotes the presheaf of categories on $\Sch$ assigning to any scheme $X$ the category of syntomic morphisms $\fromto{Y}{X}$, then the argument above also shows that nilimmersions are $\Syn$-equivalences.
\end{nul}

\begin{prp}\label{prp:redlocnoethcase} Suppose $X$ and $Y$ are reduced noetherian local schemes; then any finite universal homeomorphism $f\colon\fromto{Y}{X}$ induces an essentially surjective functor $\fromto{\Et_X}{\Et_Y}$.
\begin{proof} Write $X=\Spec A$ and $Y=\Spec B$ for reduced noetherian local rings $(A,\mathfrak{m}_A,\kappa_A)$ and $(B,\mathfrak{m}_B,\kappa_B)$. It suffices to show that any \'etale $B$-scheme $Z$ descends to an \'etale $A$-scheme.

When $\dim X=0$, the ring $A=\kappa_A$ is a field, $B=\kappa_B$ is a finite purely inseparable field extension of $A$, and $Z\cong\Spec\left(\Prod_{i\in I}F_i\right)$ corresponds to a finite product of finite separable field extensions $F_i$ of $\kappa_B$. For any $i\in I$, the subfield $E_i\subset F_i$ of elements separable over $\kappa_A$ gives a defines a finite separable field extension of $\kappa_A$ such that $F_i$ is purely inseparable over $E_i$. Each tensor product $E_i\otimes_{\kappa_A}\kappa_B$ is isomorphic to $F_i$, and the proof is complete in this case.

We now proceed by noetherian induction. Suppose now that the result is known for any $X$ of dimension $\dim X<n$. One may write $Z=Z_0\sqcup Z_1$, where $Z_0$ is finite over $X$, and the fiber of $Z_1$ over the unique closed point of $X$ is empty. Since $Z_1$ is \'etale over the complement over the closed point, the induction hypothesis applies. One is thus reduced to the case in which $Z$ is finite \'etale over $B$. In this case, $Z=\Spec D$ for some finite \'etale $B$-algebra $D$.

Now consider the completions $\mathfrak{A}$ of $A$ along $\mathfrak{m}_A$, $\mathfrak{B}$ of $B$ along $\mathfrak{m}_B$, and $\mathfrak{D}$ of $D$ along $\mathfrak{m}_BD$. It follows from \cite[18.3.4]{MR39:220} along with the dimension zero case above that $\mathfrak{D}$ descends to a finite \'etale $\mathfrak{A}$-algebra $\mathfrak{C}$. Now by \cite[5.4.4]{MR36:177c}, there is a finite $A$-algebra $C$ (essentially unique by \cite[5.4.1]{MR36:177c}) whose completion along $\mathfrak{m}_AC$ is $\mathfrak{C}$. The $A$-algebra $C$ is flat, since both $\fromto{A}{\mathfrak{A}}$ and $\fromto{C}{\mathfrak{C}}$ are faithfully flat; hence by \cite[18.4.14]{MR39:220}, $C$ is a finite \'etale $A$-algebra.
\end{proof}
\end{prp}

\begin{thm}\label{thm:UHsareEtequivs} Any universal homeomorphism is an $\Et$-equivalence.
\begin{proof} Suppose $f\colon\fromto{Y}{X}$ a universal homeomorphism. The claim is that $f$ is an $\Et$-equivalence. We employ the earlier results of this section to make a series of reductions to \ref{prp:redlocnoethcase}. (1) By \ref{cor:etequivloconbas}, one may assume that $X=\Spec A$ is affine, in which case $Y=\Spec B$ is as well. (2) Applying \ref{prp:uhlimcofiltptype} and \ref{cor:limitarguhsareets}, one may assume that $f$ is of finite type.  (3) Applying \ref{prp:factorfiniteuh} and \ref{prp:nilEtsurj}, it suffices to assume that $f$ is of finite presentation. (4) Writing $A$ as a filtered colimit of finitely generated $\ZZ$-algebras, we obtain a description of $X$ as a cofiltered limit of schemes of finite type over $\Spec\ZZ$, and employing \ref{thm:fpresuhsdescend} along with \cite[8.8.2(ii)]{MR36:178} and \cite[17.7.8]{MR39:220}, one is thus reduced to the case in which $A$ is noetherian. (5) By \ref{cor:checkonlocals}, one is further reduced to the case in which $A$ is local and noetherian, in which case $B$ is also. (6) By \ref{prp:nilEtsurj}, one may assume that both $A$ and $B$ are reduced. (7) In view of \ref{prp:uhinduceeqetale}, it suffices to show that $\fromto{\Et_X}{\Et_Y}$ is essentially surjective. This is \ref{prp:redlocnoethcase}.
\end{proof}
\end{thm}

\begin{cor} A morphism is a universal homeomorphism if and only if it is a universal $\eet$-equivalence.
\begin{proof} This is \ref{thm:UHsareEtequivs} combined with \ref{lem:Etequivsarehomeos}.
\end{proof}
\end{cor}

\section{Universal homeomorphisms and Frobenii}\label{sect:uhFrob}

\begin{ntn} Suppose $X$ an $\FF\!_p$-scheme, and suppose $Y$ an $X$-scheme. Let $q=p^r$. Denote by $\phi_X^{(q)}\colon\fromto{X}{X}$ the absolute $q$-th power Frobenius \cite[XIV=XV \S 1]{MR0491704}, and by $Y^{(q)}_{\!\!/X}$ the pullback
\begin{equation*}
\begin{tikzpicture} 
\matrix(m)[matrix of math nodes, 
row sep=3ex, column sep=5ex, 
text height=1.5ex, text depth=0.25ex] 
{Y^{(q)}_{\!\!/X}&Y\\ 
X&X\\}; 
\path[>=stealth,->,font=\scriptsize] 
(m-1-1) edge (m-1-2) 
edge (m-2-1) 
(m-1-2) edge (m-2-2) 
(m-2-1) edge node[below]{$\phi_X^{(q)}$} (m-2-2); 
\end{tikzpicture}
\end{equation*}
The absolute Frobenius $\phi_Y^{(q)}\colon\fromto{Y}{Y}$ thus factors through the \emph{relative} or \emph{geometric} Frobenius
\begin{equation*}
\Phi_{Y/X}^{(q)}\colon\fromto{Y}{Y^{(q)}_{\!\!/X}}.
\end{equation*}
\end{ntn}

\begin{dfn} An $\FF\!_p$-scheme $X$ is said to be \emph{perfect} if the absolute $p$-th power Frobenius $\phi_X^{(p)}\colon\fromto{X}{X}$ is an isomorphism. Denote by $(\Sch_{\perf}/\FF\!_p)$ the full subcategory of $(\Sch/\FF\!_p)$ comprised of perfect $\FF\!_p$-schemes.
\end{dfn}

\begin{prp} The inclusion $\into{(\Sch_{\perf}/\FF\!_p)}{(\Sch/\FF\!_p)}$ admits a right adjoint.
\begin{proof} Let $X_{\perf}$ be the limit of the diagram
\begin{equation*}
\cdots\ \tikz[baseline]\draw[>=stealth,->,font=\scriptsize](0,0.5ex)--node[above]{$\phi_X$}(0.75,0.5ex);\ X\ \tikz[baseline]\draw[>=stealth,->,font=\scriptsize](0,0.5ex)--node[above]{$\phi_X$}(0.75,0.5ex);\ X.\qedhere
\end{equation*}
\end{proof}
\end{prp}

\begin{prp}[\protect{\cite[XIV=XV \S 1 n{\textordmasculine} 2, Pr. 2(a)]{MR0491704}}] The geometric (and hence the absolute) $q$-th power Frobenius is a universal homeomorphism.
\end{prp}

\begin{cor}[\protect{\cite[XIV=XV \S 1 n{\textordmasculine} 2, Pr. 2(c)]{MR0491704}}]\label{prp:tfaefrob} Suppose $X$ and $Y$ two $\FF\!_p$-schemes. The following are equivalent for a morphism $f\colon\fromto{Y}{X}$ locally of finite presentation.
\begin{enumerate}[(\ref{prp:tfaefrob}.1)]
\item The morphism $f$ is unramified (respectively, \'etale).
\item The geometric Frobenius $\Phi^{(p)}_{Y/X}$ is unramified (resp., \'etale).
\item The geometric Frobenius $\Phi^{(p)}_{Y/X}$ is a monomorphism (resp., isomorphism).
\end{enumerate}
\end{cor}

\begin{cor}\label{cor:perfectetlfp} A morphism of perfect $\FF\!_p$-schemes is \'etale if and only if it is locally of finite presentation.
\begin{proof} This follows from \ref{prp:tfaefrob}, since the geometric Frobenius is an isomorphism.
\end{proof}
\end{cor}

\begin{prp}[(Koll{\'a}r) \protect{\cite[6.6]{MR1432036}}]\label{prp:kollar} Suppose $X$ a quasicompact $\FF_{\!p}$-scheme, and suppose $f\colon\fromto{V}{U}$ a finite universal homeomorphism of noetherian $X$-schemes. Then there is a $q=p^r$ and a morphism $f^{(-q)}\colon\fromto{U}{V^{(q)}_{\!\!/X}}$ such that the following diagram commutes:
\begin{equation*}
\begin{tikzpicture} 
\matrix(m)[matrix of math nodes, 
row sep=4ex, column sep=4ex, 
text height=1.5ex, text depth=0.25ex] 
{&U&\\
V&&V^{(q)}_{\!\!/X}.\\}; 
\path[>=stealth,->,font=\scriptsize] 
(m-2-1) edge node[above left,inner sep=1pt]{$f$} (m-1-2) 
(m-1-2) edge node[above right,inner sep=1pt]{$f^{(-q)}$} (m-2-3) 
(m-2-1) edge node[below]{$\phi_X^{(q)}$} (m-2-3); 
\end{tikzpicture}
\end{equation*}
Moreover, the morphism $f^{(-q)}$ is functorial in the sense that for any commutative square
\begin{equation*}
\begin{tikzpicture} 
\matrix(m)[matrix of math nodes, 
row sep=4ex, column sep=4ex, 
text height=1.5ex, text depth=0.25ex] 
{V&U\\ 
T&S\\}; 
\path[>=stealth,->,font=\scriptsize] 
(m-1-1) edge node[above]{$f$} (m-1-2) 
edge (m-2-1) 
(m-1-2) edge (m-2-2) 
(m-2-1) edge node[below]{$g$} (m-2-2); 
\end{tikzpicture}
\end{equation*}
of noetherian $X$-schemes in which the horizontal morphisms $f$ and $g$ are universal homeomorphisms, there is a $q=p^r$ and a commutative diagram
\begin{equation*}
\begin{tikzpicture} 
\matrix(m)[matrix of math nodes, 
row sep=0.25ex, column sep=6ex, 
text height=1.5ex, text depth=0.25ex] 
{V&&V^{(q)}_{\!\!/X}\\
&U&\\
[4ex]&S&\\
T&&T^{(q)}_{\!\!/X}.\\}; 
\path[>=stealth,->,font=\scriptsize] 
(m-1-1) edge node[above]{$\Phi^{(q)}_{V/X}$} (m-1-3) 
edge node[below,sloped,inner sep=1pt]{$f$} (m-2-2)
edge (m-4-1)
(m-2-2) edge node[below,sloped,inner sep=1pt]{$f^{(-q)}$} (m-1-3)
edge (m-3-2)
(m-3-2) edge node[above,sloped,inner sep=1pt]{$g^{(-q)}$}(m-4-3)
(m-1-3) edge (m-4-3) 
(m-4-1) edge node[below]{$\Phi^{(q)}_{T/X}$} (m-4-3)
edge node[above,sloped,inner sep=1pt]{$g$} (m-3-2); 
\end{tikzpicture}
\end{equation*}
\begin{proof}[Sketch of proof] The functoriality of $f^{(-q)}$ and the quasicompactness of $X$ permit us to assume $X$, $U$, and $V$ affine; write $X=\Spec C$, $U=\Spec A$, and $V=\Spec B$. By \ref{lem:uhschdomnilimm} it suffices to assume that $f^{\sharp}$ is either surjective or injective; i.e., that $f\colon\fromto{V}{U}$ is either a nilimmersion or schematically dominant.

Suppose first that $\fromto{V}{U}$ is a nilimmersion; let $J$ be the kernel of $f^{\sharp}$. Since $A$ is noetherian, $J$ is nilpotent; if $v=p^u$ is sufficiently large, then $J^v=0$. For such $q$ and any element $b\in B$, set $(f^{(-v)})^{\sharp}(b)=a^v$ for some preimage $a$ of $b$ (it does not matter which). This defines the desired morphism $f^{(-v)}\colon\fromto{U}{V^{(v)}_{\!\!/X}}$.

Now suppose that $\fromto{V}{U}$ is schematically dominant. We construct $f^{(-q)}$ inductively. If $U$ is an Artin scheme, then the result is immediate. One thus finds a $q=p^r$ such that $A\subset B_{/C}^{(q)}A$ is an isomorphism on generic points. (Here $B_{/C}^{(q)}$ is the $C$-algebra generated by the image of the absolute $q$-th power Frobenius.) Let $I$ denote the conductor of this extension; then by induction there is a $t=p^s$ such that $(B_{/C}^{(q)}A/I)_{/C}^{(t)}\subset(B/I)$. It follows that $B_{/C}^{(qt)}\subset(B_{/C}^{(q)}A/I)_{/C}^{(t)}\subset(B/I)$, whence one defines the desired morphism $f^{(-qt)}\colon\fromto{U}{V^{(qt)}_{\!\!/X}}$.
\end{proof}
\end{prp}

\section{Topological rigidity}\label{sect:toprigidity}

\setcounter{nul}{-1}

\begin{ntn} For any scheme $X$, denote by $(\Sch_{\coh}/X)$ the category of coherent $X$-schemes.
\end{ntn}

\begin{dfn}\label{dfn:toprigid} A coherent scheme $Y$ is said to be \emph{topologically rigid} if any universal homeomorphism $\fromto{Y'}{Y}$ with $Y'$ reduced is an isomorphism. Denote by $(\Sch_{\trig}/X)$ the full subcategory of $(\Sch_{\coh}/X)$ spanned by topologically rigid $X$-schemes; denote by $j\colon\into{(\Sch_{\trig}/X)}{(\Sch_{\coh}/X)}$ the inclusion.
\end{dfn}

\begin{lem} A scheme $Y$ is topologically rigid if and only if any universal homeomorphism $\fromto{Y'}{Y}$ admits a section.
\begin{proof} If $Y$ is topologically rigid, then for any universal homeomorphism $\fromto{Y'}{Y}$, the composite $\fromto{Y'_{\red}}{Y}$ is an isomorphism. Conversely, a section of a universal homeomorphism $\fromto{Y'}{Y}$ is a nilimmersion, hence an isomorphism if $Y'$ is reduced.
\end{proof}
\end{lem}

\begin{lem} A topologically rigid scheme is reduced.
\begin{proof} If $Y$ is topologically rigid, then $\fromto{Y_{\red}}{Y}$ admits a section.
\end{proof}
\end{lem}

\begin{lem}\label{lem:coprodtoprig} The coproduct of a set $\{U_i\}_{i\in I}$ of schemes is topologically rigid if and only if each scheme $U_i$ is so.
\begin{proof} This is immediate from the observation that any coproduct of universal homeomorphisms is a universal homeomorphism.
\end{proof}
\end{lem}

\begin{lem}\label{lem:trigsuffffp} A scheme $Z$ is topologically rigid if and only if $Z$ is reduced, and, for any reduced scheme $V$, any universal homeomorphism of finite presentation $\fromto{V}{Z}$ admits a section.
\begin{proof} The necessity is obvious. For the sufficiency, suppose $Z$ satisfies the two conditions, and suppose $f\colon\fromto{W}{Z}$ a universal homeomorphism, with $W$ reduced; the claim is that $f$ is an isomorphism. Apply \ref{prp:uhlimcofiltptype} to write $W=\lim_{\alpha\in\Lambda^{\op}}W_{\alpha}$ for a cofiltered diagram of coherent schemes over $Z$ under $W$, where each morphism $\fromto{W_{\alpha}}{Z}$ is finite. One may \cite[8.7.1]{MR36:178} assume that for each $\alpha\in M$, the scheme $W_{\alpha}$ is reduced. Now apply \ref{prp:factorfiniteuh} for each $\alpha\in\Lambda$ to factor $\fromto{W_{\alpha}}{Z}$ as a nilimmersion $\fromto{W_{\alpha}}{W_{\alpha}'}$ followed by a universal homeomorphism $\fromto{W_{\alpha}'}{Z}$ of finite presentation. Now $\fromto{W_{\alpha}'}{Z}$ admits a section by assumption, and since $Z$ is reduced, this section factors through $W_{\alpha}$, whence $\fromto{W_{\alpha}}{Z}$ is an isomorphism. Thus $f$ is a limit of isomorphisms, verifying the claim.
\end{proof}
\end{lem}

\begin{prp}\label{prp:jcommuteswithfiltlims} Suppose $X$ a coherent scheme, and suppose $\Lambda$ a filtered category. Then the inclusion $j\colon\into{(\Sch_{\trig}/X)}{(\Sch_{\coh}/X)}$ reflects limits of $\Lambda^{\op}$-diagrams of affine $X$-schemes.
\begin{proof} Suppose $Z$ a coherent $X$-scheme with the property that $Z=\lim_{\alpha\in\Lambda^{\op}}Z_{\alpha}$, where each $Z_{\alpha}$ is a topologically rigid affine $X$-scheme. We employ \ref{lem:trigsuffffp}: suppose now that $W$ is a reduced scheme, and that $f\colon\fromto{W}{Z}$ is a universal homeomorphism of finite presentation. By \ref{thm:fpresuhsdescend} and \cite[8.8.2(ii)]{MR36:178}, there exists an element $\alpha\in\Lambda$ such that $\fromto{W}{Z}$ descends to a universal homeomorphism $\fromto{W_{\alpha}}{Z_{\alpha}}$. This universal homeomorphism admits a section, since $Z_{\alpha}$ is topologically rigid, and so $\fromto{W}{Z}$ admits a section as well.
\end{proof}
\end{prp}

\begin{ctrexm} Observe that $j$ does not commute with finite limits, for even the tensor product of perfect fields may have nilpotents.
\end{ctrexm}

\begin{nul} Topological rigidity is related to weak normality [\S \ref{sect:wknorm}], \cite{MR556308,MR811809,MR1432036,MR1998022}, and for reduced and irreducible affines, it is equivalent to absolute weak normality as in \cite[\S B]{rydh:submers}. The related geometric notion is that of \emph{seminormality}, which coincides with topological rigidity in characteristic zero. In effect, a scheme is seminormal if and only if all non-normalities are maximally transverse.
\end{nul}

\begin{dfn} A universal homeomorphism $f\colon\fromto{Y'}{Y}$ of schemes is said to be \emph{arithmetically trivial}\footnote{Traverso et al. \cite{MR0277542,MR571055} use the word ``quasi-isomorphism'' for such universal homeomorphisms.} if for any point $y'\in Y'$, if $y=f(y')$, then the field extension $\into{\kappa(y)}{\kappa(y')}$ is an isomorphism. A coherent scheme $Y$ is said to be \emph{seminormal} \cite{MR595029} if any arithmetically trivial universal homeomorphism $\fromto{Y'}{Y}$ with $Y'$ reduced is an isomorphism.
\end{dfn}

\begin{lem} Any seminormal scheme is reduced.
\end{lem}

\begin{exm}\label{exm:twnsch} We leave the proofs of these facts as exercises for the amused reader; they will not be used in the sequel.
\begin{enumerate}[(\ref{exm:twnsch}.1)]
\item Suppose $k$ a field. Then $\Spec k$ is always seminormal, but it is topologically rigid if and only if $k$ is perfect.
\item Suppose $X$ a normal scheme. Then $X$ is seminormal as well. If, for any generic point $\eta$ of $X$, the field $\mathscr{O}_{X,\eta}=\kappa(\eta)$ is perfect, then $X$ is topologically rigid.
\item\label{item:snequiv} Suppose $X$ a reduced scheme whose set of irreducible components is locally finite. Let $\mathscr{R}_{X}$ denote the quasicoherent $\mathscr{O}_X$-algebra of rational functions on $X$ \cite[8.3.3]{EGA1}. Then the following are equivalent (cf. \cite[1.4 and 1.7]{MR618807}).
\begin{enumerate}[(\ref{exm:twnsch}.\ref{item:snequiv}.1)]
\item The scheme $X$ is seminormal.
\item For every point $x\in X$, the local scheme $\Spec\mathscr{O}_{X,x}$ is seminormal.
\item For every point $x\in X$, and for any integral element $f\in\mathscr{R}_{X,x}$ over $\mathscr{O}_{X,x}$, the conductor of $\mathscr{O}_{X,x}$ in $\mathscr{O}_{X,x}[f]$ is a radical ideal of $\mathscr{O}_{X,x}[f]$.
\item For every point $x\in X$, if $f\in\mathscr{R}_{X,x}$ is an integral element with the property that both $f^u,f^v\in\mathscr{O}_{X,x}$ for positive relatively prime natural numbers $u$ and $v$, then $f\in\mathscr{O}_{X,x}$.
\item For every point $x\in X$, if $f\in\mathscr{R}_{X,x}$ is an integral element with the property that both $f^2,f^3\in\mathscr{O}_{X,x}$, then $f\in\mathscr{O}_{X,x}$.
\end{enumerate}
\item\label{item:reirredwnorm} Suppose $X$ a reduced scheme whose set of irreducible components is locally finite; then $X$ is topologically rigid if and only $X$ satisfies the following conditions \cite[Th. 1]{MR714467}.
\begin{enumerate}[(\ref{exm:twnsch}.\ref{item:reirredwnorm}.1)]
\item The scheme $X$ is seminormal.
\item For every point $x\in X$ and every integral element $f\in\mathscr{R}_{X,x}^{\perf}$ over $\mathscr{O}_{X,x}$ with the property that there is a prime number $p$ such that $f^p\in\mathscr{O}_{X,x}$ and $pf\in\mathscr{O}_{X,x}$, one has $f\in\mathscr{O}_{X,x}$.
\end{enumerate}
\item Suppose $X$ a reduced scheme whose set of irreducible components is locally finite. Then if $X$ is geometrically unibranch, it is seminormal if and only if it is normal.
\item Suppose $k$ a field. Then a nodal curve (e.g., $y^2=x^2$ or $y^2=x^3+x^2$) over $k$ is seminormal.
\item Suppose $k$ a field. Then for any nonnegative integer $n$, the union of the $n$ coordinate axes in $\mathbf{A}_k^{\!n}$ is a seminormal scheme.
\end{enumerate}
\end{exm}

\begin{ctrexm}\label{ctrexm:twnsch} In effect, a scheme can fail to be topologically rigid for two reasons. First, there is the failure of seminormality.
\begin{enumerate}[(\ref{ctrexm:twnsch}.1)]
\item Suppose $k$ a field. Neither cuspidal curves (e.g., $y^2=x^3$) nor tacnodal curves (e.g., $y^2=x^4$) nor ramphoid cusps (e.g., $y^2=x^5$), etc. ..., are seminormal.
\item Suppose $k$ a field. Then for any nonnegative integer $n$, the scheme comprised of $m$ lines intersecting in one point in $\mathbf{A}_k^{\!n}$ is not seminormal if $m>n$.
\suspend{enumerate}
Beyond the ``geometric'' considerations of seminormality, there is also the ``arithmetic'' requirement that the function field of each irreducible component be perfect.
\resume{enumerate}[{[(\ref{ctrexm:twnsch}.1)]}]
\item Suppose $k$ a field of characteristic $p$; suppose $n$ a positive integer. Then consider the subring
\begin{equation*}
A\coloneq k[x^{p^n},x^jy\ :\ 0\leq j<p^n]\subset k[x,y].
\end{equation*}
Since the induced morphism $\fromto{\mathbf{A}_k^{\!2}}{\Spec A}$ is a nontrivial universal homeomorphism, it follows immediately that $\Spec A$ is not topologically rigid.
\item An $\FF_p$-scheme is topologically rigid only if it is perfect. In particular [\ref{cor:perfectetlfp}], if $X$ is of finite type over a perfect field $k$, then $X$ is topologically rigid only if it is \emph{\'etale} over $\Spec k$.
\item If $k$ is a perfect field of characteristic $p$, then no smooth $k$-scheme of positive dimension is topologically rigid, but any of them is seminormal.
\end{enumerate}
\end{ctrexm}

\section{Topological rigidification}\label{sect:trig}

\begin{prp}\label{prp:trigexists} Suppose $X$ a coherent scheme. The inclusion functor
\begin{equation*}
j\colon\into{(\Sch_{\trig}/X)}{(\Sch_{\coh}/X)}
\end{equation*}
admits a right adjoint.
\begin{proof} For any $X$-scheme $Y$, denote by $(\UH/Y)$ the full subcategory of $(\Sch/Y)$ comprised of $Y$-schemes $Z$ such that the structural morphism $\fromto{Z}{Y}$ is a universal homeomorphism. This category contains an essentially ($\UU$-)small full cofinal subcategory comprised of the finite universal homeomorphisms $\fromto{Z}{Y}$ [\ref{prp:uhlimcofiltptype}]. Let $Y_{\trig}$ denote the limit of the diagram $\fromto{(\UH/Y)}{(\Sch/Y)}$. This is a cofiltered diagram of affine $Y$-schemes since finite limits exist in $(\UH/Y)$, so \ref{prp:cofiltlimuhuh} provides a natural universal homeomorphism $\fromto{Y_{\trig}}{Y}$. 

The $X$-scheme $Y_{\trig}$ is topologically rigid: indeed, suppose $f\colon\fromto{Z}{Y_{\trig}}$ a universal homeomorphism; then the composite $\fromto{Z}{Y}$ is a universal homeomorphism, so the projection $\sigma\colon\fromto{Y_{\trig}}{Z}$ is a section of $f$.

Now define a functor in the following manner: assign to any $X$-scheme $Y$ the object $Y_{\trig}$ of $(\UH/Y)$. For any morphism $g\colon\fromto{Y'}{Y}$ of $X$-schemes, let
\begin{equation*}
g_{\trig}\colon\fromto{Y'_{\trig}}{Y_{\trig}}
\end{equation*}
be the composite on the left:
\begin{equation*}
\begin{tikzpicture} 
\matrix(m)[matrix of math nodes, 
row sep=4ex, column sep=4ex, 
text height=1.5ex, text depth=0.25ex] 
{Y'_{\trig}&Y_{\trig}\times_YY'&Y_{\trig}\\
&Y'&Y.\\}; 
\path[>=stealth,->,font=\scriptsize] 
(m-1-1) edge (m-1-2)
edge (m-2-2)
(m-1-2) edge (m-2-2) 
edge (m-1-3) 
(m-2-2) edge node[below]{$g$} (m-2-3) 
(m-1-3) edge (m-2-3); 
\end{tikzpicture}
\end{equation*}

If $Y$ is an $X$-scheme, there is a canonical universal homeomorphism $\epsilon\colon\fromto{Y_{\trig}}{Y}$; this will be the counit of the adjunction. If $Y'$ is a topologically rigid $X$-scheme, then the natural morphism $\fromto{Y'_{\trig}}{Y'}$ is an isomorphism; the inverse $\eta$ will be the unit of the adjunction. We now show that the unit and counit give an isomorphism
\begin{equation*}
\Mor_X(Y',Y_{\trig})\cong\Mor_X(Y',Y),
\end{equation*}
thereby completing the proof.

Suppose now $Y$ and $Y'$ two $X$-schemes; suppose $Y'$ topologically rigid. If $g\colon\fromto{Y'}{Y}$ a morphism of $X$-schemes, one sees easily that the composite
\begin{equation*}
Y'\ \tikz[baseline]\draw[>=stealth,->,font=\scriptsize](0,0.5ex)--node[above]{$\eta$}(0.5,0.5ex);\ Y'_{\trig}\ \tikz[baseline]\draw[>=stealth,->](0,0.5ex)--(0.5,0.5ex);\ Y_{\trig}\ \tikz[baseline]\draw[>=stealth,->,font=\scriptsize](0,0.5ex)--node[above]{$\epsilon$}(0.5,0.5ex);\ Y
\end{equation*}
is again $g$. On the other hand, if $f\colon\fromto{Y'}{Y_{\trig}}$ is a morphism of $X$-schemes, then one forms the corresponding map by $\fromto{Y'}{Y}$ composing with $\epsilon$, and one returns to a morphism $\fromto{Y'}{Y_{\trig}}$ by forming the following diagram:
\begin{equation*}
\begin{tikzpicture} 
\matrix(m)[matrix of math nodes, 
row sep=4ex, column sep=4ex, 
text height=1.5ex, text depth=0.25ex] 
{Y'_{\trig}&Y_{\trig}\times_YY'&[4ex]Y_{\trig}\times_YY_{\trig}&Y_{\trig}\\
&Y'&Y_{\trig}&Y.\\}; 
\path[>=stealth,->,font=\scriptsize] 
(m-1-1) edge (m-1-2)
edge node[below left,inner sep=1pt]{$\eta^{-1}$} (m-2-2)
(m-1-2) edge (m-2-2) 
edge node[above]{$\id\times f$} (m-1-3) 
(m-2-2) edge node[below]{$f$} (m-2-3) 
(m-1-3) edge (m-2-3)
(m-1-3) edge (m-1-4) 
(m-2-3) edge node[below]{$\epsilon$} (m-2-4) 
(m-1-4) edge node[right]{$\epsilon$} (m-2-4); 
\end{tikzpicture}
\end{equation*}
In order to show that $\eta$ composed with the long composite on the top is equal to $f$, it suffices to show that the square
\begin{equation*}
\begin{tikzpicture} 
\matrix(m)[matrix of math nodes, 
row sep=4ex, column sep=5ex, 
text height=1.5ex, text depth=0.25ex] 
{Y'_{\trig}&Y_{\trig}\times_YY'\\
Y_{\trig}&Y_{\trig}\times_YY_{\trig}\\}; 
\path[>=stealth,->,font=\scriptsize] 
(m-1-1) edge (m-1-2) 
edge (m-2-1) 
(m-1-2) edge node[right]{$\id\times f$} (m-2-2) 
(m-2-1) edge node[below]{$\Delta$} (m-2-2); 
\end{tikzpicture}
\end{equation*}
commutes. This is immediate, as $\Delta$ is a nilimmersion, and $Y'_{\trig}$ is reduced.
\end{proof}
\end{prp}

\begin{dfn} We call the adjoint $(-)_{\trig}$ to $j$ the \emph{topological rigidification}.
\end{dfn}

\begin{lem} The topological rigidification of any coherent scheme $Y$ is the limit of the natural diagram of $Y$-schemes indexed by $(\UH_{\red}/Y):=(\UH/Y)\cap(\Sch_{\red}/Y)$.
\begin{proof} The functor $\goesto{Z}{Z_{\red}}$ induces a right adjoint to the inclusion $\into{(\UH_{\red}/Y)}{(\UH/Y)}$; hence $(\UH_{\red}/Y)$ is cofinal in $(\UH/Y)$.
\end{proof}
\end{lem}

\begin{thm}\label{exm:trigvsperf} Suppose $p$ a prime number. For any reduced, noetherian $\FF_p$-scheme $X$, the topological rigidification of $X$ is isomorphic to the perfection $X_{\perf}$.
\begin{proof} Since $\fromto{X_{\perf}}{X}$ is a limit of universal homeomorphisms, it is a universal homeomorphism [\ref{prp:cofiltlimuhuh}]. We employ \ref{lem:trigsuffffp}: suppose $W$ a reduced scheme, and suppose $g\colon\fromto{W}{X_{\perf}}$ a finite universal homeomorphism of finite presentation. We wish to construct a section of $g$. By \ref{thm:fpresuhsdescend} and \cite[8.8.2(ii)]{MR36:178}, $W$ descends to a universal homeomorphism $f\colon\fromto{V}{X}$ of finite presentation. Now apply \ref{prp:kollar}: there exists a $q=p^r$ and a morphism $f^{(-q)}\colon\fromto{U}{V^{(q)}_{\!\!/X}}$ such that the following diagram commutes:
\begin{equation*}
\begin{tikzpicture} 
\matrix(m)[matrix of math nodes, 
row sep=4ex, column sep=4ex, 
text height=1.5ex, text depth=0.25ex] 
{&X&\\
V&&V^{(q)}_{\!\!/X}.\\}; 
\path[>=stealth,->,font=\scriptsize] 
(m-2-1) edge node[above left,inner sep=1pt]{$f$} (m-1-2) 
(m-1-2) edge node[above right,inner sep=1pt]{$f^{(-q)}$} (m-2-3) 
(m-2-1) edge node[below]{$\phi_X^{(q)}$} (m-2-3); 
\end{tikzpicture}
\end{equation*}
Now form the pullback of $f$ along the Frobenius $\phi_X^{(q)}$:
\begin{equation*}
\begin{tikzpicture} 
\matrix(m)[matrix of math nodes, 
row sep=3ex, column sep=5ex, 
text height=1.5ex, text depth=0.25ex] 
{V^{(q)}_{\!\!/X}&V\\ 
X&X.\\}; 
\path[>=stealth,->,font=\scriptsize] 
(m-1-1) edge (m-1-2) 
edge node[left]{$f^{(q)}$} (m-2-1) 
(m-1-2) edge node[right]{$f$} (m-2-2) 
(m-2-1) edge node[below]{$\phi_X^{(q)}$} (m-2-2); 
\end{tikzpicture}
\end{equation*}
The identities
\begin{equation*}
f^{(q)}\circ\Phi^{(q)}_{/X}=f\textrm{\quad and\quad}f^{(-q)}\circ f=\Phi^{(q)}_{/X}
\end{equation*}
together imply that $f^{(q)}\circ f^{(-q)}\circ f=f$, but since $f$ is an epimorphism in the category of schemes [\ref{prp:schdomuhsareepis}], one deduces that $f^{(-q)}$ is a section of $f^{(q)}$. Since $g$ is the pullback of $f^{(q)}$, it admits a section as well.
\end{proof}
\end{thm}

\section{Weak normalization}\label{sect:wknorm}

\begin{dfn} For any morphism $f\colon\fromto{Y}{X}$ of schemes, let $\UH_f$ denote the category of factorizations
\begin{equation*}
Y\ \tikz[baseline]\draw[>=stealth,->](0,0.5ex)--(0.5,0.5ex);\ X'\ \tikz[baseline]\draw[>=stealth,->](0,0.5ex)--(0.5,0.5ex);\ X
\end{equation*}
of $f$ with the property that $\fromto{X'}{X}$ is a universal homeomorphism. If $\UH_f$ consists of the terminal object alone, then we will say that $X$ is \emph{weakly normal under $f$} or \emph{weakly normal under $Y$}. The initial object of $\UH_f$, if it exists, will be called the \emph{weak normalization $X_{\wn,f}$ of $X$ under $Y$} (or the \emph{weak normalization relative to $f$} if there is possibility for confusion). Thus a scheme $X$ is weakly normal under a morphism $f\colon\fromto{Y}{X}$ if and only if the weak normalization exists and coincides with $X$.
\end{dfn}

\begin{prp} Suppose $f\colon\fromto{Y}{X}$ a morphism of coherent schemes. Suppose $\fromto{X'}{X}$ a faithfully flat morphism, and set $Y'\coloneq X'\times_XY$. If $X'$ is weakly normal under $Y'$, then $X$ is weakly normal under $Y$.
\begin{proof} Suppose $X'$ weakly normal under $Y'$. Suppose
\begin{equation*}
Y\ \tikz[baseline]\draw[>=stealth,->](0,0.5ex)--(0.5,0.5ex);\ Z\ \tikz[baseline]\draw[>=stealth,->](0,0.5ex)--(0.5,0.5ex);\ X
\end{equation*}
a factorization of $f$ with the property that $\fromto{Z}{X}$ is a universal homeomorphism. Then the projection $\fromto{Z'\coloneq Z\times_XX'}{X'}$ is an isomorphism. But by faithfully flat descent, it follows that $\fromto{Z}{X}$ is an isomorphism as well.
\end{proof}
\end{prp}

\begin{thm}[(Andreotti--Bombieri) \protect{\cite[Teorema 2]{MR0266923}}] Suppose $X$ a reduced, coherent scheme. Suppose $f\colon\fromto{Y}{X}$ a dominant (hence schematically dominant \cite[5.4.3]{EGA1}), coherent morphism of schemes. Then the weak normalization of $X$ under $Y$ exists.
\begin{proof}[Sketch of proof] Write $\mathscr{A}\coloneq f_{\star}\mathscr{O}_Y$, which is a quasicoherent $\mathscr{O}_X$-algebra. Since the affination $\Spec_X\mathscr{A}$ of $Y$ over $X$ is dominant, and since universal homeomorphisms are affine, one may assume that $Y$ is affine over $X$.

For any point $x\in X$, let $p=\chr\kappa(x)$, and set
\begin{equation*}
{}^{\mathrm{wn},f}\mathscr{O}_{X,x}\coloneq\{a\in\mathscr{A}_x\ |\ \exists m\in\NN,\ a^{p^m}\in f^{\sharp}_x(\mathscr{O}_{X,x})+\mathscr{J}(\mathscr{A}_x)\},
\end{equation*}
where $\mathscr{J}$ denotes the Jacobson radical. Now define, for any open set $U$ of $X$,
\begin{equation*}
{}^{\mathrm{wn},f}\mathscr{O}_X(U)\coloneq\{a\in\mathscr{A}(U)\ |\ \forall x\in U,\ a_x\in{}^{\mathrm{wn},f}\mathscr{O}_{X,x}\}.
\end{equation*}
One verifies that this defines a quasicoherent $\mathscr{O}_X$-subalgebra ${}^{\mathrm{wn},f}\mathscr{O}_X$ of $\mathscr{A}$ \cite[Proposizione 4]{MR0266923}.

Now set $X_{\wn,f}\coloneq\Spec_X({}^{\mathrm{wn},f}\mathscr{O}_X)$. The first claim is that the morphism ${}^{\wn}f\colon\fromto{X_{\wn,f}}{X}$ is a universal homeomorphism. Indeed, it is a straightforward matter to verify that ${}^{\mathrm{wn},f}\mathscr{O}_X$ is an integral $\mathscr{O}_X$-algebra and that ${}^{\wn}f$ is surjective. To show that ${}^{\wn}f$ is universally injective, in effect, one shows that, for any point $x\in X$, the ring ${}^{\mathrm{wn},f}\mathscr{O}_{X,x}$ is a local ring whose maximal ideal is the Jacobson radical $\mathscr{J}(\mathscr{A}_x)$. To verify that ${}^{\wn}f$ is universally injective at $x$, one wishes to show that the field extension
\begin{equation*}
\fromto{\kappa(x)}{({}^{\mathrm{wn},f}\mathscr{O}_{X,x}/\mathscr{J}(\mathscr{A}_x))}
\end{equation*}
is purely inseparable. This follows from the observation that an element $a\in({}^{\mathrm{wn},f}\mathscr{O}_{X,x}/\mathscr{J}(\mathscr{A}_x))$ only if for some nonnegative integer $m$, one has $a^{p^m}\in\kappa(x)$.

Now suppose $\mathscr{B}$ a quasicoherent $\mathscr{O}_X$-subalgebra of $\mathscr{A}$ with the property that the morphism $\fromto{\Spec_X\mathscr{B}}{X}$ is a universal homeomorphism. Then for any point $x\in X$, it is easy to see that $\mathscr{B}_x$ is contained in the $\mathscr{O}_{X,x}$-subalgebra of $\mathscr{A}_x$ generated by $\mathscr{J}(\mathscr{A}_x)$. From this it follows that $X_{\wn,f}$ has the desired universal property.
\end{proof}
\end{thm}

\begin{prp}\label{prp:trigiswn} Suppose $\fromto{Y}{X}$ a dominant morphism of reduced, coherent schemes. Then the weak normalization of $X$ under the composite $Y_{\trig}\ \tikz[baseline]\draw[>=stealth,->](0,0.5ex)--(0.5,0.5ex);\ Y\ \tikz[baseline]\draw[>=stealth,->](0,0.5ex)--(0.5,0.5ex);\ X$ coincides with the topological rigidification of $X$.
\begin{proof} Write $f$ for the composite $Y_{\trig}\ \tikz[baseline]\draw[>=stealth,->](0,0.5ex)--(0.5,0.5ex);\ Y\ \tikz[baseline]\draw[>=stealth,->](0,0.5ex)--(0.5,0.5ex);\ X$. Then $\fromto{X_{\wn,f}}{X}$ is a universal homeomorphism, and for any universal homeomorphism $g\colon\fromto{X'}{X_{\wn,f}}$, the morphism $\fromto{Y_{\trig}}{Y}$ factors through $X'\times_XY$ and thus through $X'$ itself, whence $g$ admits a section. Hence $X_{\wn,f}$ is the topological rigidification.
\end{proof}
\end{prp}

\begin{ntn*} For the next result, we use the following notations. Suppose $X$ a reduced scheme whose set of irreducible components is locally finite. Let $X^{(0)}$ be the set of generic points of irreducible components of $X$, and set
\begin{equation*}
\Xi:=\Coprod_{\eta\in X^{(0)}}\Spec\left(\kappa(\eta)^{\perf}\right);
\end{equation*}
then the natural morphism $\pi\colon\fromto{\Xi}{X}$ is dominant.
\end{ntn*}

\begin{cor}\label{cor:computetrig} The weak normalization of $X$ under $\pi\colon\fromto{\Xi}{X}$ is the topological rigidification of $X$.
\begin{proof} Apply \ref{prp:trigiswn}, noting that by \ref{lem:coprodtoprig}, $\Xi$ is topologically rigid.
\end{proof}
\end{cor}

\begin{nul} Suppose $X$ a reduced scheme whose set of irreducible components is finite. Then one may form ``the'' \emph{absolute integral closure} $\overline{X}$ of $X$ \cite{MR0289501,MR0224600,MR0257064}. Then one can show that the topological rigidification of $X$ is isomorphic to its weak normalization under $\fromto{\overline{X}}{X}$. This construction of $X_{\trig}$ is less satisfying, however, since it only recovers the noncanonical isomorphism type.
\end{nul}

\begin{prp}[(Manaresi), \protect{\cite[I.6]{MR556308}}]\label{prp:manaresi} Suppose $A$ and $A'$ reduced rings, and suppose $\into{A}{A'}$ an integral extension. Then the weak normalization of $\Spec A$ under $\Spec A'$ is $\Spec A''$, where $A''$ is the equalizer
\begin{equation*}
A''\ \tikz[baseline]\draw[>=stealth,->](0,0.5ex)--(0.5,0.5ex);\ A'\ \begin{tikzpicture}[baseline] \draw[>=stealth,->] (0,1ex) -- (0.5,1ex); \draw[>=stealth,->] (0,0.25ex) -- (0.5,0.25ex); \end{tikzpicture}\ (A'\otimes_AA')_{\red}
\end{equation*}
in the category of rings.
\begin{proof}[Sketch of proof] Write ${}^{A'/\wn}\!A$ for the ring such that $\Spec({}^{A'/\wn}\!A)$ is weak normalization of $\Spec A$ under $\Spec A'$. One verifies directly that $\fromto{\Spec A''}{\Spec A}$ is radicial and surjective, hence a universal homeomorphism, so we have the inclusion $A''\subset{}^{A'/\wn}\!A$. On the other hand, every element $a\in{}^{A'/\wn}\!A$ has the property that $a\otimes 1-1\otimes a$ is a nilpotent element of $A'\otimes_AA'$, so we also have the inclusion $A''\supset{}^{A'/\wn}\!A$.
\end{proof}
\end{prp}

\begin{cor}\label{lem:etalelocwn} Suppose
\begin{equation*}
\begin{tikzpicture} 
\matrix(m)[matrix of math nodes, 
row sep=4ex, column sep=4ex, 
text height=1.5ex, text depth=0.25ex] 
{Y'&X'\\ 
Y&X\\}; 
\path[>=stealth,->,font=\scriptsize] 
(m-1-1) edge node[above]{$h'$} (m-1-2) 
edge node[left]{$g$} (m-2-1) 
(m-1-2) edge node[right]{$f$} (m-2-2) 
(m-2-1) edge node[below]{$h$} (m-2-2); 
\end{tikzpicture}
\end{equation*}
a pullback square of reduced coherent schemes. Suppose also that $f$ and $g$ are dominant, and that $h$ is \'etale. Then the square
\begin{equation*}
\begin{tikzpicture} 
\matrix(m)[matrix of math nodes, 
row sep=4ex, column sep=5ex, 
text height=1.5ex, text depth=0.25ex] 
{Y_{g/\wn}&X_{f/\wn}\\ 
Y&X\\}; 
\path[>=stealth,->,font=\scriptsize] 
(m-1-1) edge (m-1-2) 
edge (m-2-1) 
(m-1-2) edge (m-2-2) 
(m-2-1) edge node[below]{$h$} (m-2-2); 
\end{tikzpicture}
\end{equation*}
is a pullback square.
\begin{proof} The question is local, so assume $X$ and $Y$ affine. Form the normalizations $X_{f/\mathrm{n}}$ and $Y_{g/\mathrm{n}}$ of $X$ and $Y$ under $W$; by \cite[18.12.15]{MR39:220}, one has $Y_{g/\mathrm{n}}\cong X_{f/\mathrm{n}}\times_XY$, whence we may assume $f$ and $g$ integral. Write $X=\Spec A$, $Y=\Spec B$, $X'=\Spec A'$, and $Y'=\Spec B'$. Hence $X_{f/\wn}=\Spec A''$ for $A''$ the equalizer in \ref{prp:manaresi}. By the flatness of $\into{A}{B}$ combined with \cite[17.5.7]{MR39:220}, the diagram
\begin{equation*}
A''\otimes_AB\ \tikz[baseline]\draw[>=stealth,->](0,0.5ex)--(0.5,0.5ex);\ B'\ \begin{tikzpicture}[baseline] \draw[>=stealth,->] (0,1ex) -- (0.5,1ex); \draw[>=stealth,->] (0,0.25ex) -- (0.5,0.25ex); \end{tikzpicture}\ (B'\otimes_{A'}B')_{\red}
\end{equation*}
is an equalizer as well.
\end{proof}
\end{cor}

\begin{cor}\label{cor:etcovtrigistrig} Suppose $X$ and $Y$ reduced, coherent schemes. Suppose $p\colon\fromto{Y}{X}$ an \'etale surjective morphism. Then the following square is a pullback square
\begin{equation*}
\begin{tikzpicture} 
\matrix(m)[matrix of math nodes, 
row sep=4ex, column sep=4ex, 
text height=1.5ex, text depth=0.25ex] 
{Y_{\trig}&Y\\ 
X_{\trig}&X.\\}; 
\path[>=stealth,->,font=\scriptsize] 
(m-1-1) edge (m-1-2) 
edge node[left]{$p_{\trig}$} (m-2-1) 
(m-1-2) edge node[right]{$p$} (m-2-2) 
(m-2-1) edge (m-2-2); 
\end{tikzpicture}
\end{equation*}
\begin{proof} By absolute noetherian approximation \cite[C.9]{MR92f:19001} combined with \ref{prp:jcommuteswithfiltlims}, \cite[8.8.2]{MR36:178}, \cite[17.7.8]{MR39:220}, and \cite[8.7.1]{MR36:178}, one reduces to the case in which $X$ is of finite type over $\Spec\ZZ$. By \ref{prp:trigiswn} and \ref{lem:etalelocwn}, it is enough to find a topologically rigid scheme $W$ and a dominant morphism $\fromto{W}{X}$ such that $W\times_XY$ is topologically rigid and $\fromto{W\times_XY}{Y}$ is dominant. For this, let $X^{(0)}$ denote the finite set of generic points of the irreducible components of $X$, and set
\begin{equation*}
W:=\Coprod_{\eta\in X^{(0)}}\Spec\left(\Omega(\eta)\right)
\end{equation*}
for some choice of algebraic closure $\Omega(\eta)$ of $\kappa(\eta)$ for each $\eta\in X^{(0)}$.
\end{proof}
\end{cor}

\begin{ctrexm} Suppose $X$, $Y$, and $Y'$ reduced coherent schemes; suppose $\fromto{Y}{X}$ an \'etale morphism; and suppose $\fromto{Y'}{Y}$ a universal homeomorphism. In general there is not a universal homeomorphism $\fromto{X'}{X}$ and an identification of $Y'\cong X'\times_XY$. To see this, let $X=\Spec k$ for some imperfect field $k$, $Y=\Spec (k\times k)$, and $Y'=\Spec (k^{\perf}\times k)$.
\end{ctrexm}

\section{The h site and the étale site of topologically rigid schemes}\label{sect:htopos}

\setcounter{nul}{-1}

\begin{nul} In this section, the language of $\infty$-categories and $\infty$-topoi \cite{lurie_inftytopoi} will be used systematically. In particular, we do not distinguish notationally an ordinary category from its nerve.
\end{nul}

\begin{ntn} Suppose $X$ a coherent scheme. The inclusion
\begin{equation*}
j\colon\into{(\Sch_{\trig}/X)}{(\Sch_{\coh}/X)}
\end{equation*}
induces a string of adjoints (left to right) $(j_{\natural},j^{\star},j_{\star},j^{\natural})$ between the $\infty$-categories of $\VV$-small presheaves of simplicial sets on $(\Sch_{\coh}/X)$ and on $(\Sch_{\trig}/X)$, where for any presheaf $E$ of simplicial sets on $(\Sch_{\coh}/X)$ and any presheaf $\Theta$ on $(\Sch_{\trig}/X)$, one has
\begin{equation*}
j_{\star}E=E((-)_{\trig})\qquad\textrm{and}\qquad j^{\star}\Theta=\Theta(j(-)),
\end{equation*}
$j_{\natural}E$ is the usual left Kan extension of $F$ along $j$, and $j^{\natural}\Theta$ is the right Kan extension of $G$ along $\goesto{Y}{Y_{\trig}}$.\footnote{These notations are compatible with those of Giraud \cite[Chap. 0, \S 3]{MR49:8992} if one defines $j^{-1}$ as the functor $\goesto{Y}{Y_{\trig}}$.}
\end{ntn}

\begin{dfn} The $\infty$-topos of sheaves of $\VV$-small simplicial sets \cite[6.2.2.6]{lurie_inftytopoi} on the category $(\Sch_{\coh}/X)$ equipped with the \'etale topology will be called the \emph{large \'etale $\infty$-topos} and will be denoted
\begin{equation*}
(\EEt/X)\coloneq\Shv_{\et}(\Sch_{\coh}/X).
\end{equation*}
(This is not a $\VV$-small category; it is rather a topos in the universe $\VV$ \cite[6.3.5.17]{lurie_inftytopoi}.)

The \'etale topology can be restricted to $(\Sch_{\trig}/X)$; the resulting $\infty$-topos of sheaves of $\VV$-small simplicial sets \cite[6.2.2.6]{lurie_inftytopoi} will be called the \emph{topologically invariant étale $\infty$-topos} (also an $\infty$-topos in $\VV$) and will be denoted
\begin{equation*}
(\EEt_{\trig}/X)\coloneq\Shv_{\et}(\Sch_{\trig}/X).
\end{equation*}
\end{dfn}

\begin{lem}\label{lem:keyV} Suppose $X$ a coherent scheme. Then the functor $\goesto{Y}{Y_{\trig}}$ induces a morphism of sites
\begin{equation*}
j\colon\fromto{(\Sch_{\trig}/X)_{\et}}{(\Sch_{\coh}/X)_{\et}}
\end{equation*}
and thus a geometric morphism of $\infty$-topoi
\begin{equation*}
j\colon\fromto{(\EEt_{\trig}/X)}{(\EEt/X)}.
\end{equation*}
\begin{proof} By \ref{cor:etcovtrigistrig} if $\{\fromto{V_i}{U}\}_{i\in I}$ is an \'etale covering of an $X$-scheme $U$, then each morphism $\fromto{V_{i,\trig}}{U_{\trig}}$ is \'etale, and since $\fromto{U_{\trig}}{U}$ is a universal homeomorphism, $\fromto{\Coprod_{i\in I}V_{i,\trig}}{U_{\trig}}$ is faithfully flat. Thus $j$ is continuous.

For any \'etale sheaf $H$ over $X$, the presheaf $j^{\star}H$ is already an h sheaf. Hence $j^{\star}$ admits a further left adjoint $j_!$, which is the composite of the functor $j_{\natural}$ with \'etale sheafification. In particular, $j^{\star}$ is left exact.
\end{proof}
\end{lem}

\begin{dfn} For any coherent scheme $X$, we will say that a sheaf $E\in(\EEt/X)$ is \emph{topologically invariant} if, for any universal homeomoprhism $\fromto{Y'}{Y}$ over $X$, the induced map $\fromto{E(Y)}{E(Y')}$ is an equivalence.
\end{dfn}

\begin{prp}\label{prp:hsubtopos} Suppose $X$ a coherent scheme. Then the geometric morphism $j$ is an embedding of topoi that identifies the topologically invariant étale $\infty$-topos with the full subcategory of $(\EEt/X)$ spanned by the topologically invariant sheaves.
\begin{proof} It is immediate that $j^{\star}j_{\star}\simeq\id$, so $j$ is indeed an embedding of $\infty$-topoi; moreover, for any \'etale sheaf $E$, the natural morphism $\fromto{E}{j_{\star}j^{\star}E}$ is an equivalence if and only if, for any $X$-scheme $Y$, the induced map $\fromto{E(Y)}{E(Y_{\trig})}$ is an equivalence.
\end{proof}
\end{prp}

\begin{cor} For any morphism  $f\colon\fromto{X}{Y}$ of coherent schemes, there is a morphism of topoi $f\colon\fromto{(\EEt_{\trig}/X)}{(\EEt_{\trig}/Y)}$ such that the following diagrams commute up to homotopy:
\begin{equation*}
\begin{tikzpicture}[baseline]
\matrix(m)[matrix of math nodes, 
row sep=6ex, column sep=5ex, 
text height=1.5ex, text depth=0.25ex] 
{(\EEt_{\trig}/X)&(\EEt_{\trig}/Y)\\
(\EEt/X)&(\EEt/Y)\\}; 
\path[>=stealth,->,font=\scriptsize] 
(m-1-1) edge node[above]{$f_{\star}$} (m-1-2) 
edge[right hook->] node[left]{$j_{\star}$} (m-2-1) 
(m-1-2) edge[left hook->] node[right]{$j_{\star}$} (m-2-2) 
(m-2-1) edge node[below]{$f_{\star}$} (m-2-2); 
\end{tikzpicture}
\textrm{\quad and\quad}
\begin{tikzpicture}[baseline]
\matrix(m)[matrix of math nodes, 
row sep=6ex, column sep=5ex, 
text height=1.5ex, text depth=0.25ex] 
{(\EEt_{\trig}/Y)&(\EEt_{\trig}/X)\\
(\EEt/Y)&(\EEt/X).\\}; 
\path[>=stealth,->,font=\scriptsize] 
(m-1-1) edge node[above]{$f^{\star}$} (m-1-2) 
edge[right hook->] node[left]{$j_{\star}$} (m-2-1) 
(m-1-2) edge[left hook->] node[right]{$j_{\star}$} (m-2-2) 
(m-2-1) edge node[below]{$f^{\star}$} (m-2-2); 
\end{tikzpicture}
\end{equation*}
\begin{proof} Using \ref{prp:hsubtopos} to regard the topologically invariant étale $\infty$-topos as a full subcategory of the \'etale $\infty$-topos, the claims are: (1) that for any topologically invariant sheaf $E\in(\EEt/X)$, the sheaf $f_{\star}E$ is topologically invariant, and (2) that for any topologically invariant sheaf $H\in(\EEt/Y)$, the sheaf $f^{\star}E$ is topologically invariant. The first claim follows from the fact that pullbacks of universal homeomorphisms are universal homeomorphisms, and the second is obvious.
\end{proof}
\end{cor}

\begin{lem}\label{lem:jstarYisYtrig} Suppose $X$ a coherent scheme. The inverse image functor \begin{equation*}
j^{\star}\colon\fromto{(\EEt/X)}{(\EEt_{\trig}/X)}
\end{equation*}
sends representables to representables; in particular, one has an isomorphism of h sheaves $i^{\star}Y\cong Y_{\trig}$ for any coherent $X$-scheme $Y$.
\begin{proof} Suppose $Y$ an $X$-scheme. Then for any h sheaf $E$, one has
\begin{equation*}
\Map_{\EEt_{\trig}}(Y_{\trig},E)\simeq E(Y_{\trig})=j_{\star}E(Y)\simeq\Map_{\EEt}(Y,j_{\star}E)\simeq\Map_{\EEt_{\trig}}(j^{\star}Y,E),
\end{equation*}
so the result follows from Yoneda.
\end{proof}
\end{lem}

\begin{nul} We have constructed a string of adjoints $(j_!,j^{\star},j_{\star})$ relating the \'etale and the topologically invariant étale $\infty$-topoi. Perhaps contrary to expectation, the functor $j_!$ does not provide an identification of the h $\infty$-topos $(\EEt_{\trig}/X)$ with a slice $\infty$-topos $((\EEt/X)/Y)$; that is, $j$ is not an \'etale morphism of topoi in the sense of \cite[6.3.5]{lurie_inftytopoi}. Nevertheless, the functor $j_!$ has a clear geometric meaning, since it preserves representables; in particular, for any \'etale sheaf $E$, one may write $H_{\trig}$ for $j_!j^{\star}E$, since the endofunctor $j_!j^{\star}$ is the left Kan extension of $\goesto{Y}{Y_{\trig}}$ along the Yoneda embedding. This permits one to say when a large \'etale sheaf (e.g., a Deligne--Mumford $n$-stack) is topologically rigid.
\end{nul}

\begin{dfn} A large \'etale sheaf $H\in(\EEt/X)$ over a scheme $X$ is said to be \emph{topologically rigid} if the natural morphism $\fromto{E_{\trig}}{E}$ is an equivalence.
\end{dfn}

\begin{prp} For any coherent scheme $X$, the functor $j_!\colon\fromto{(\EEt_{\trig}/X)}{(\EEt/X)}$ is a fully faithful functor that identifies $(\EEt_{\trig}/X)$ with the full subcategory of $(\EEt/X)$ spanned by the topologically rigid sheaves.
\begin{proof} The claim is that the counit $\fromto{\id}{j^{\star}j_!}$ is an equivalence. Since both $j^{\star}$ and $j_!$ commute with colimits, it is enough to show that for any topologically rigid $X$-scheme $Y$, the morphism $\fromto{Y}{j^{\star}j_!Y}$ is an equivalence. Since the \'etale topology is subcanonical, $j_!Y$ is isomorphic to $jY$, and the result follows from \ref{lem:jstarYisYtrig}. 
\end{proof}
\end{prp}

\begin{dfn} The topology induced on $(\Sch_{\coh}/X)$ by the topologically invariant \'etale $\infty$-topos \cite[6.2.4.2]{lurie_inftytopoi} is called the \emph{topologically invariant \'etale topology}.
\end{dfn}

\begin{nul} The h topology of Voevodsky and Suslin--Voevodsky \cite{MR1403354,MR1376246} is in fact strictly finer than the topologically invariant \'etale topology, but the sheafifications of representables with respect to these two topologies coincide. To formulate this, fix an excellent scheme $X$, and denote by $(\Sch^{\sim,\trig}/X)$ the full subcategory of the $\infty$-topos $(\EEt_{\trig}/X)$ spanned by those sheaves that are represented by the topological rigidification of a scheme of finite type over $X$, and denote by $(\Sch^{\sim,\mathrm{h}}/X)$ the full subcategory of the $\infty$-topos $\Shv_{\mathrm{h}}(\Sch_{\mathrm{ft}}/X)$ spanned by those sheaves obtained as the sheafification of representable presheaves.
\end{nul}

\begin{thm} The $\infty$-categories $(\Sch^{\sim,\trig}/X)$ and $\Shv_{\mathrm{h}}(\Sch_{\mathrm{ft}}/X)$ are canonically equivalent.
\begin{proof} This is now a formal consequence of \cite[3.2.9]{MR1403354}.
\end{proof}
\end{thm}

\begin{cor} Topologically rigid schemes over $X$ form a generating site for the h $\infty$-topos for which the topology is subcanonical.
\end{cor}

\begin{cor}[cf. Suslin--Voevodsky, \protect{\cite[\S 10]{MR1376246}}]\label{cor:hcohomology} Suppose $G$ a group scheme over $X$. Then $G_{\trig}$ is also a group scheme over $X$, and for any scheme $Y$ of finite type over $X$, there are canonical isomorphisms of (nonabelian) cohomology
\begin{equation*}
H^{\ast}_{\mathrm{h}}(Y,G)\cong H^{\ast}_{\et}(Y,G_{\trig}).
\end{equation*}
\begin{proof} Since $(-)_{\trig}$ is a right adjoint, it commutes with fiber products, so $G_{\trig}$ is a group. The mapping space from $Y$ to $G$ in the h $\infty$-topos is homotopy equivalent to the mapping space from $Y_{\trig}$ to $A_{\trig}$ in the topologically invariant $\infty$-topos, which in turn coincides with the \'etale cohomology of $Y$ with $A$ coefficients by the \frenchquote{invariance topologique} of the \'etale topos [\ref{thm:UHsareEtequivs}].
\end{proof}
\end{cor}

\bibliographystyle{amsalpha}
\bibliography{../../math,../../egasga}

\newcommand{\noopsrt}[1]{}
\providecommand{\bysame}{\leavevmode\hbox to3em{\hrulefill}\thinspace}
\providecommand{\MR}{\relax\ifhmode\unskip\space\fi MR }
\providecommand{\MRhref}[2]{%
  \href{http://www.ams.org/mathscinet-getitem?mr=#1}{#2}
}
\providecommand{\href}[2]{#2}
\begin{thebibliography}{{SGA}64b}

\bibitem[AB69]{MR0266923}
A.~Andreotti and E.~Bombieri, \emph{Sugli omeomorfismi delle variet\`a
  algebriche}, Ann. Scuola Norm. Sup Pisa (3) \textbf{23} (1969), 431--450.
  \MR{MR0266923 (42 \#1825)}

\bibitem[Art71]{MR0289501}
M.~Artin, \emph{On the joins of {H}ensel rings}, Advances in Math. \textbf{7}
  (1971), 282--296 (1971). \MR{MR0289501 (44 \#6690)}

\bibitem[Con07]{MR2356346}
Brian Conrad, \emph{Deligne's notes on {N}agata compactifications}, J.
  Ramanujan Math. Soc. \textbf{22} (2007), no.~3, 205--257. \MR{MR2356346
  (2009d:14002)}

\bibitem[EGA$_{\textrm{I}}^{\ast}$]{EGA1}
A.~Grothendieck and J.~A. Dieudonn\'{e}, \emph{El\'{e}ments de
  g\'{e}om\'{e}trie alg\'{e}brique {I}}, Grundlehren der math. Wissenschaften
  in Einzeldarstellungen, vol. 166, Springer-Verlag, Heidelberg, 1971.

\bibitem[EGA$_{\textrm{II}}$]{MR36:177b}
A.~Grothendieck, \emph{\'{E}l\'ements de g\'eom\'etrie alg\'ebrique. {II}.
  \'{E}tude globale \'el\'ementaire de quelques classes de morphismes}, Inst.
  Hautes \'Etudes Sci. Publ. Math. (1961), no.~8, 222. \MR{36 \#177b}

\bibitem[EGA$_{\textrm{III}}$ I]{MR36:177c}
\bysame, \emph{\'{E}l\'ements de g\'eom\'etrie alg\'ebrique. {III}. \'{E}tude
  cohomologique des faisceaux coh\'erents. {I}}, Inst. Hautes \'Etudes Sci.
  Publ. Math. (1961), no.~11, 167. \MR{36 \#177c}

\bibitem[EGA$_{\textrm{IV}}$ II]{MR33:7330}
\bysame, \emph{\'{E}l\'ements de g\'eom\'etrie alg\'ebrique. {IV}. \'{E}tude
  locale des sch\'emas et des morphismes de sch\'emas. {II}}, Inst. Hautes
  \'Etudes Sci. Publ. Math. (1965), no.~24, 231. \MR{33 \#7330}

\bibitem[EGA$_{\textrm{IV}}$ III]{MR36:178}
\bysame, \emph{\'{E}l\'ements de g\'eom\'etrie alg\'ebrique. {IV}. \'{E}tude
  locale des sch\'emas et des morphismes de sch\'emas. {III}}, Inst. Hautes
  \'Etudes Sci. Publ. Math. (1966), no.~28, 255. \MR{36 \#178}

\bibitem[EGA$_{\textrm{IV}}$ IV]{MR39:220}
\bysame, \emph{\'{E}l\'ements de g\'eom\'etrie alg\'ebrique. {IV}. \'{E}tude
  locale des sch\'emas et des morphismes de sch\'emas {IV}}, Inst. Hautes
  \'Etudes Sci. Publ. Math. (1967), no.~32, 361. \MR{39 \#220}

\bibitem[Eno68]{MR0224600}
Edgar Enochs, \emph{Totally integrally closed rings}, Proc. Amer. Math. Soc.
  \textbf{19} (1968), 701--706. \MR{MR0224600 (37 \#199)}

\bibitem[Gir71]{MR49:8992}
Jean Giraud, \emph{Cohomologie non ab\'elienne}, Die Grundlehren der
  mathematischen Wissenschaften, Band 179, Springer-Verlag, Berlin, 1971.
  \MR{49 \#8992}

\bibitem[Gro97a]{MR1483108}
Alexander Grothendieck, \emph{Brief an {G}. {F}altings}, Geometric {G}alois
  actions, 1, London Math. Soc. Lecture Note Ser., vol. 242, Cambridge Univ.
  Press, Cambridge, 1997, With an English translation on pp. 285--293,
  pp.~49--58. \MR{MR1483108 (99c:14023)}

\bibitem[Gro97b]{MR1483107}
Alexandre Grothendieck, \emph{Esquisse d'un programme}, Geometric {G}alois
  actions, 1, London Math. Soc. Lecture Note Ser., vol. 242, Cambridge Univ.
  Press, Cambridge, 1997, With an English translation on pp. 243--283,
  pp.~5--48. \MR{MR1483107 (99c:14034)}

\bibitem[GT80]{MR571055}
S.~Greco and C.~Traverso, \emph{On seminormal schemes}, Compositio Math.
  \textbf{40} (1980), no.~3, 325--365. \MR{MR571055 (81j:14030)}

\bibitem[Hoc70]{MR0257064}
M.~Hochster, \emph{Totally integrally closed rings and extremal spaces},
  Pacific J. Math. \textbf{32} (1970), 767--779. \MR{MR0257064 (41 \#1718)}
  
\bibitem[HTT]{lurie_inftytopoi}
J.~Lurie, \emph{Higher topos theory}, Annals of Mathematics Studies, no. 170,
  Princeton University Press, 2009, Also available electronically from
  \textit{http://www.math.harvard.edu/$\sim$lurie}.

\bibitem[Ill71]{MR58:10886a}
L.~Illusie, \emph{Complexe cotangent et d\'eformations. {I}}, Springer-Verlag,
  Berlin, 1971, Lecture Notes in Mathematics, Vol. 239. \MR{58 \#10886a}

\bibitem[Kol97]{MR1432036}
J{\'a}nos Koll{\'a}r, \emph{Quotient spaces modulo algebraic groups}, Ann. of
  Math. (2) \textbf{145} (1997), no.~1, 33--79. \MR{MR1432036 (97m:14013)}

\bibitem[LV81]{MR618807}
John~V. Leahy and Marie~A. Vitulli, \emph{Seminormal rings and weakly normal
  varieties}, Nagoya Math. J. \textbf{82} (1981), 27--56. \MR{MR618807
  (83a:14015)}

\bibitem[Man80]{MR556308}
Mirella Manaresi, \emph{Some properties of weakly normal varieties}, Nagoya
  Math. J. \textbf{77} (1980), 61--74. \MR{MR556308 (81e:13011)}

\bibitem[Pic03]{MR1998022}
Gabriel Picavet, \emph{Universally going-down rings, 1-split rings, and
  absolute integral closure}, Comm. Algebra \textbf{31} (2003), no.~10,
  4655--4685. \MR{MR1998022 (2004d:13010)}

\bibitem[Ryd07]{rydh:submers}
D.~Rydh, \emph{Submersions and effective descent of {\'e}tale morphisms}, To
  appear in Bull. Soc. Math. France. Also available from
  \textit{arXiv:math/0710.2888v3}, 2007.

\bibitem[SGA$_1$]{MR50:7129}
\emph{Rev\^etements \'etales et groupe fondamental}, S\'eminaire de
  G\'eom\'etrie Alg\'ebrique du Bois Marie 1960--61 (SGA 1). Dirig\'e par A.
  Grothendieck. Lecture Notes in Mathematics, Vol. 224, Springer-Verlag,
  Berlin, 1960--61. \MR{50 \#7129}

\bibitem[SGA$_4$ I]{MR50:7130}
\emph{Th\'eorie des topos et cohomologie \'etale des sch\'emas. {T}ome 1:
  {T}h\'eorie des topos}, S\'eminaire de G\'eom\'etrie Alg\'ebrique du Bois
  Marie 1963--64 (SGA 4). Dirig\'e par M. Artin, A. Grothendieck, J. L.
  Verdier. Avec la collaboration de N. Bourbaki, P. Deligne, B. Saint--Donat.
  Lecture Notes in Mathematics, Vol. 269, Springer-Verlag, Berlin, 1963--64.
  \MR{50 \#7130}

\bibitem[SGA$_4$ II]{MR50:7131}
\emph{Th\'eorie des topos et cohomologie \'etale des sch\'emas. {T}ome 2},
  S\'eminaire de G\'eom\'etrie Alg\'ebrique du Bois Marie 1963--64 (SGA 4).
  Dirig\'e par M. Artin, A. Grothendieck, J. L. Verdier. Avec la collaboration
  de N. Bourbaki, P. Deligne, B. Saint--Donat. Lecture Notes in Mathematics,
  Vol. 270, Springer-Verlag, Berlin, 1963--64. \MR{50 \#7131}

\bibitem[SGA$_5$]{MR0491704}
\emph{Cohomologie {$\ell$}-adique et fonctions {$L$}}, S\'eminaire de
  G\'eom\'etrie Alg\'ebrique du Bois Marie 1965--66 (SGA 5). Dirig\'e par A.
  Grothendieck. Avec la collaboration de I. Bucur, C. Houzel, L. Illusie,
  J.--P. Jouanolou, et J.--P. Serre. Lecture Notes in Mathematics, Vol. 589,
  Springer-Verlag, Berlin, 1965--66. \MR{MR0491704 (58 \#10907)}

\bibitem[Sti02]{MR1931963}
Jakob Stix, \emph{Affine anabelian curves in positive characteristic},
  Compositio Math. \textbf{134} (2002), no.~1, 75--85. \MR{MR1931963
  (2003i:14027)}

\bibitem[SV96]{MR1376246}
Andrei Suslin and Vladimir Voevodsky, \emph{Singular homology of abstract
  algebraic varieties}, Invent. Math. \textbf{123} (1996), no.~1, 61--94.
  \MR{MR1376246 (97e:14030)}

\bibitem[Swa80]{MR595029}
Richard~G. Swan, \emph{On seminormality}, J. Algebra \textbf{67} (1980), no.~1,
  210--229. \MR{MR595029 (82d:13006)}

\bibitem[Tra70]{MR0277542}
Carlo Traverso, \emph{Seminormality and {P}icard group}, Ann. Scuola Norm. Sup.
  Pisa (3) \textbf{24} (1970), 585--595. \MR{MR0277542 (43 \#3275)}

\bibitem[TT90]{MR92f:19001}
R.~W. Thomason and Thomas Trobaugh, \emph{Higher algebraic {$K$}-theory of
  schemes and of derived categories}, The Grothendieck Festschrift, Vol.\ III,
  Progr. Math., vol.~88, Birkh\"auser Boston, Boston, MA, 1990, pp.~247--435.
  \MR{92f:19001}

\bibitem[Voe90]{MR1098621}
V.~A. Voevodski{\u\i}, \emph{\'{E}tale topologies of schemes over fields of
  finite type over {${\bf Q}$}}, Izv. Akad. Nauk SSSR Ser. Mat. \textbf{54}
  (1990), no.~6, 1155--1167. \MR{MR1098621 (92j:14023)}

\bibitem[Voe96]{MR1403354}
V.~Voevodsky, \emph{Homology of schemes}, Selecta Math. (N.S.) \textbf{2}
  (1996), no.~1, 111--153. \MR{MR1403354 (98c:14016)}

\bibitem[Yan83]{MR714467}
Hiroshi Yanagihara, \emph{Some results on weakly normal ring extensions}, J.
  Math. Soc. Japan \textbf{35} (1983), no.~4, 649--661. \MR{MR714467
  (85e:13010)}

\bibitem[Yan85]{MR811809}
\bysame, \emph{On an intrinsic definition of weakly normal rings}, Kobe J.
  Math. \textbf{2} (1985), no.~1, 89--98. \MR{MR811809 (87d:13007)}

\end{thebibliography}

\end{document}